\documentclass[12pt]{article}
\usepackage{graphicx,amsthm}
\usepackage{latexsym,amsfonts}
\newcommand{\eqdef}{\, =\kern -12.7pt\raise 6pt\hbox{{\tiny\textrm{def}}}\,\,}
\newcommand{\bgs}{\\\smallskip \\}
\newtheorem{theorem}{Theorem}
\newtheorem{lemma}{Lemma}
\newtheorem{corollary}{Corollary}

\def\hypf#1#2#3#4#5{\, _{#1}\kern -1pt F_{#2}\kern -3pt\left[
\begin{matrix}#3\\#4\end{matrix}\,;\,#5\right]}
\title{How to lose as little as possible}
\author{Vittorio Addona, Stan Wagon, and Herb Wilf}
\date{}
\bigskip
\begin{document}
\maketitle

\begin{abstract}
 Suppose Alice has a coin with heads probability $q$ and Bob has one with heads probability $p>q$.
  Now each of them will toss their coin $n$ times, and Alice will win iff she gets more heads than Bob does. Evidently the game favors Bob, but for the given $p,q$, what is the choice of $n$ that maximizes Alice's  chances of winning? We show that there is an essentially unique value $N(q,p)$ of $n$ that maximizes the probability $f(n)$ that the weak coin will win, and it satisfies $\left\lfloor{\frac{1}{2(p-q)}-\frac12}\right\rfloor\le N(q,p)\le \left\lceil{\frac{\max{(1-p,q)}}{p-q}}\right\rceil$. The analysis uses the multivariate form of Zeilberger's algorithm  to find an indicator function $J_n(q,p)$ such that $J>0$ iff $n<N(q,p)$ followed by a close study of this function, which is a linear combination of two Legendre polynomials. An integration-based  algorithm is given for computing $N(q,p)$.
\end{abstract}

\newpage
{\small\parskip 0pt
\tableofcontents }
\newpage

\section{The problem}
 Suppose Alice has a coin with heads probability $q$ and Bob has one with heads probability $p$.
 Suppose $q<p$.  Now each of them will toss their coin $n$ times, and Alice wins iff she gets more heads
  than Bob does (n.b.: in case of a tie, Bob wins). Evidently the game favors Bob, but for the given $p,q$, what is the choice of $n$ that
   maximizes Alice's  chances of winning?

\begin{figure}[htp]\label{fig:fig1}
\begin{center}
\includegraphics{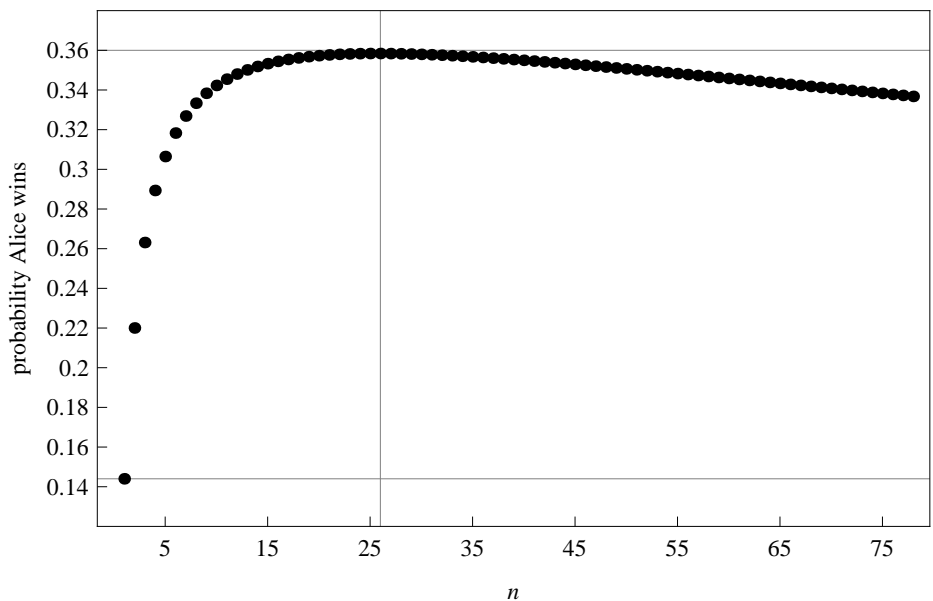}
\end{center}
\caption{Probability that Alice wins vs. $n$.}
\end{figure}

  Interestingly, there is a nontrivial (i.e., in general $>1$) unique value of $n$ that  maximizes her probability of winning. For example, in the case $p=0.2$, $q=0.18$, Figure \ref{fig:fig1} is a plot of Alice's win probability as a function of $n$.

\bigskip

In this example, if each player flips their coin 26 times, which is the best choice for her, Alice's chance of winning will be about $0.36$, compared to a chance of $0.14$ if each coin is tossed only once.

In general, her chances of winning are
\begin{equation}
\label{eq:eq1}
f(n)=f(n,p,q)\eqdef\sum_{r\ge 0}{n\choose r}p^r(1-p)^{n-r}\sum_{s> r}{n\choose s}q^s(1-q)^{n-s}.
\end{equation}
This problem, which first appeared in \cite{wa}, arose from a consideration of real-world events in the National Football League, where teams play a season of 16 games and do not play all other teams. If teams A and B have probabilities $q$ and $p>q$, respectively, of winning any game and never play each other, one can wonder about the chance that A's season record will be  strictly better than that of B. That is easy to answer, but then one is led to the question of whether the season length, 16, is favorable or not to such an outcome and what the optimal choice would be. A study of a related topic, where the central issue is the chance that the underdog beats a certain point-spread, has been published by T. Lengyel \cite{tl}.

We will also give, in section \ref{sec:sec10} an algorithm that uses repeated numerical integration to compute
the optimum value $N(q, p)$. Mathematica code for various computations, graphics, and algorithms (e.g., the generation
 of graphs of $p_n$ or computation of $N(q,p)$) is available in the electronic supplement at \cite{ws}.

Much of the work here has relied on computing power, both for numerical experiments and for proofs using symbolic computation. Some sophisticated
 algorithms in Maple and Mathematica (Zeilberger's \texttt{MultiZeil}, and cylindrical
algebra reduction of polynomial systems) played  crucial roles; without them the discoveries and proofs would have been difficult, if not impossible, to find.

\subsection{Acknowledgments} Professor Bruno Salvy, of INRIA, France, has kindly supplied to us some highly refined asymptotic results, which were of great assistance in this work. We thank Rob Knapp and Tamas Lengyel for some helpful discussions.

\section{Overview of methods and results}

It develops that there is, in this problem, a nice \textit{indicator function} $J_n(q,p)$, which is simply a linear combination of two consecutive Legendre polynomials, with the property that the sign of $f(n+1)-f(n)$ is the same as the sign of $J_n(q,p)$. We will find this indicator by using the multivariate form of Zeilberger's algorithm \cite{mz}. We will then show that for small $n$, $J$ is positive and for large enough $n$, $J$ is negative, and that there is only a single integer value of $n$, or a consecutive pair $(n,n+1)$, at which the sign of $J$ changes. Thus $f$ has a unique maximum, at $n=N(q,p)$, say. Here is the precise result.
\begin{theorem}
\label{th:mainth}
With $f(n)$ defined by (\ref{eq:eq1}) we have
\begin{equation}
\label{eq:fmon}
\frac{f(n+1)-f(n)}{((1-p)(1-q))^{n+1}}=\left(y+\frac12(1+xy)\right)\phi_n(xy)-\frac12 \phi_{n+1}(xy),
\end{equation}
where $x=p/(1-p)$, $y=q/(1-q)$, and
\begin{equation}
\label{eq:pdef}
\phi_n(z)=\sum_{r=0}^n{n\choose r}^2z^r=(1-z)^nP_n\left(\frac{1+z}{1-z}\right),
\end{equation}
and $P_n(t)$ is the classical Legendre polynomial. Therefore the indicator function
\[J_n(q,p)=\left(y+\frac12(1+xy)\right)\phi_n(xy)-\frac12 \phi_{n+1}(xy)\]
has the desired properties.
\end{theorem}
We remark that, once found, the recurrence (\ref{eq:fmon}) can be proved directly, i.e., without Zeilberger's algorithm, with little difficulty.

Next, in section \ref{sec:uniq} we will prove uniqueness of and find upper and lower bounds for $N(q,p)$ by using various properties of the
Legendre polynomials and by  a close study of a function $p_n(q)$ which for each $q\in (0,n/(2n+1))$, is the unique value of $p$ for which
$J_n(q,p)=0$. The properties of the curves $p=p_n(q)$ in  the $(p,q)$ plane play crucial roles here. First, concerning uniqueness, we have
\begin{theorem}\label{th:uni}(Unimodality) Given probabilities $p>q$ with $p+q\ne 1$, there are either one, or two consecutive, values of $n$ such that
\begin{enumerate}
\item $f(n)\ge f(n-1)$, and
\item $f(n+1)\le f(n)$, and
\item at least one of the above two inequalities is strict.
\end{enumerate}
\end{theorem}
\noindent \textbf{Definition}. Given $q<p$, let $N(q,p)$ be the value of $n$ that maximizes $f(n,p,q)$. When the value is not unique, define $N$ to be the smaller of the two possible values  that yield the maximum.

It follows from Theorem \ref{th:uni} and the definition of $N$ that $N(q,p)$ is the smallest integer $n$ such that $J_n(q,p)\le 0$.

The resulting upper and lower bounds for $N(q,p)$ are given by
\begin{theorem}\label{th:bounds}
If $N(q,p)$ is the choice of $n$ that maximizes the probability that the player with the weaker coin will win (and with ties going to the lower value) we have:
\begin{enumerate}
\item $\lfloor{\frac{1}{2(p-q)} -\frac12}\rfloor\le  N(q,p)$, but if $p+q\neq 1$, then $\lfloor{\frac{1}{2(p-q)} +\frac12}\rfloor \le  N(q,p)$.
\item $N(q,p)\le \left\lceil{\frac{\max{(1-p, q)}}{p - q}}\right\rceil$.
\end{enumerate}
\end{theorem}
Section \ref{sec:sec9} contains proofs of various properties of the graphs of $p_n$, and they are used to obtain improvements to the upper and lower bounds on $N$.
\subsection{Definitions and notation}
The heads probabilities of the two coins are $p$ and $q$, with $p>q$. We write $x=p/(1-p)$, $y=q/(1-q)$, $z=xy$, $u=(1+z)/(1-z)=1+2pq/(1-p-q)$, $\rho=(1-p+q)/(1-p-q)$. Further, $P_n$ is the $n$th Legendre polynomial and $r_n=r_n(u)=P_n(u)/P_{n-1}(u)$. The $P_n$'s satisfy the well known recurrence
\begin{equation}
\label{eq:prcr}
P_{n+1}(u)=\frac{2n+1}{n+1}uP_n(u)-\frac{n}{n+1}P_{n-1}(u).
\end{equation}
If we divide through by $P_n(u)$, we obtain the \textit{ratio recurrence}
\begin{equation}
\label{eq:ratrcr}
r_{n+1}(u)=\frac{2n+1}{n+1}u-\frac{n}{(n+1)r_{n}(u)}.
\end{equation}
which will be of use in the sequel.

The \textit{indicator function} $J_n(q,p)$, which has the sign of $f(n+1)-f(n)$, is
\[J_n(q,p)=y\phi_n(z)-\psi_n(z)=\frac12\frac{(1-p-q)^n}{((1-p)(1-q))^{n+1}}((1-p+q)P_n(u)-(1-p-q)P_{n+1}(u)),\]
where $\phi_n$ is given by (\ref{eq:pdef}) and\footnote{For a quick proof of (\ref{eq:fnlpsi}), square both sides of the Pascal triangle recurrence.}
\begin{equation}
\psi_n(z)=\sum_{r=0}^n{n\choose r+1}{n\choose r}z^{r+1}
=\frac12 (\phi_{n+1}(z)-(1+z)\phi_n(z))\label{eq:fnlpsi}
\end{equation}
$T$ will denote the \textit{interior} of the triangle in the $(p,q)$ plane whose vertices are $(0,0)$, $(1,0)$, $(\frac12,\frac12)$. $T_n$ will be the \textit{open} interval $(0,n/(2n+1))$. The line $L_k(q)$ in the $(p,q)$ plane is the line $p=\frac{1}{2k+1}+q$, and the line $M_n(q)$ is
\begin{equation}\label{eq:mnq} M_n(q):\ p=\frac{1}{n+1}+\frac{n}{n+1} q.\end{equation}
\section{Finding the indicator function}
Our first task will be to find a recurrence for $f(n)$. To do this we will use the multivariate
form of Zeilberger's algorithm, \texttt{MulZeil} \cite{mz}. As usual the results that are returned by the algorithm
can be easily verified by substitution.

Remarkably, this recurrence will show that $f(n+1)-f(n)$ is simply expressible in terms of Legendre
polynomials; this will enable us to identify the values of $n$ for which $f(n+1)\ge f(n)$ and those for which $f(n+1)\le f(n)$.

In view of eq. (\ref{eq:fmon}), Alice's probability of winning increases with $n$ as long as
\begin{equation}
\label{eq:find}
y\phi_n(xy)-\psi_n(xy)=\left(y+\frac12(1+xy)\right)\phi_n(xy)-\frac12 \phi_{n+1}(xy)
\end{equation}
is positive, and decreases otherwise. We will show that for fixed $x,y$ there is a unique value of $n$ at which
this function changes its sign.

\section{Finding the recurrence for $f(n)$}
In this section we will find the recurrence that is satisfied by $f(n)$, the sum in eq. (\ref{eq:eq1}), using the
 multidimensional version of Zeilberger's algorithm. This will prove Theorem \ref{th:mainth}.
\subsection{Finding the recurrence for the summand}
With
\begin{equation}
\label{eq:subs}
x=p/(1-p),\quad y=q/(1-q),\quad  g(n)=f(n)/((1-p)^n(1-q)^n),
\end{equation}
 the definition (\ref{eq:eq1}) of $f(n)$ becomes
\[g(n)=\sum_{r\ge 0}\sum_{s> r}{n\choose r}{n\choose s}x^ry^s.\]
Let $G(n,r,s)={n\choose r}{n\choose s}x^ry^s$, be the summand. We use Zeilberger's algorithm, and his program \texttt{MulZeil} returns a recurrence
\begin{eqnarray}
\label{eq:mainrecr}
G(n+1,r,s)-(x+1)(y+1)G(n,r,s)&=&(K_r-1)(c_1(n,r,s)G(n,r,s))\nonumber\\
&&\qquad +(K_s-1)(c_2(n,r,s)G(n,r,s)),
\end{eqnarray}
where $K_r,K_s$ are forward shift operators in their subscripts, and the $c_i$ are given by
\begin{equation}
\label{eq:thecs}
c_1=c_1(n,r,s)=\frac{r(1+y)}{r-n-1};\qquad c_2=c_2(n,r,s)=\frac{s(n+1)}{(s-n-1)(n-r+1)}.
\end{equation}
This is the recurrence for the summand, and it can be quickly verified by dividing through by $G(n,r,s)$, canceling all of the factorials, and noting that the resulting polynomial identity states that $0=0$.
\subsection{Finding the recurrence for the sum}
To find the recurrence for the sum, we sum the recurrence (\ref{eq:mainrecr}) over $s> r$,
and then sum the result over $r\ge 0$. To do this we have first, for every function $\phi$ of compact support,
\begin{equation}\label{eq:cs1}
\sum_{r\ge 0}\sum_{s> r}(K_r-1)\phi(r,s)=-\sum_{s\ge 1}\phi(0,s)+\sum_{r\ge 1}\phi(r,r),
\end{equation}
and
\begin{equation}
\label{eq:cs2}
\sum_{r\ge 0}\sum_{s> r}(K_s-1)\phi(r,s)=-\sum_{r\ge 0}\phi(r,r+1).
\end{equation}
Consequently, if we sum the recurrence (\ref{eq:mainrecr}) there results
\begin{eqnarray*}
g(n+1)-(x+1)(y+1)g(n)&=&-\sum_{s\ge 1}c_1(n,0,s)G(n,0,s)+\sum_{r\ge 1}c_1(n,r,r)G(n,r,r)\\
&&\qquad -\sum_{r\ge 0}c_2(n,r,r+1)G(n,r,r+1).
\end{eqnarray*}

Next we insert the values, from (\ref{eq:thecs}),
\[c_1(n,0,s)=0;\quad c_1(n,r,r)=\frac{r(1+y)}{r-n-1};\qquad c_2(n,r,r+1)=\frac{(r+1)(n+1)}{(r-n)(n-r+1)},\]
which gives
\[
g(n+1)-(x+1)(y+1)g(n)=\sum_{r\ge 1}\frac{r(1+y)}{r-n-1}G(n,r,r)
 -\sum_{r\ge 0}\frac{(r+1)(n+1)}{(r-n)(n-r+1)}G(n,r,r+1).
\]
Now substitute the values $G(n,r,r)={n\choose r}^2(u)^r$, and $G(n,r,r+1)={n\choose r}{n\choose r+1}x^ry^{r+1}$, and simplify the result, to obtain
\begin{eqnarray*}
g(n+1)-(x+1)(y+1)g(n)&=&\sum_{r\ge 1}\frac{r(1+y)}{r-n-1}{n\choose r}^2(xy)^r\\
 &&\qquad -\sum_{r\ge 0}\frac{(r+1)(n+1)}{(r-n)(n-r+1)}{n\choose r}{n\choose r+1}x^ry^{r+1}\\
&=&-(y+1)\sum_{r\ge 0}{n\choose r+1}{n\choose r}x^{r+1}y^{r+1}+\sum_{r\ge 0}{n+1\choose r}{n\choose r}x^ry^{r+1}\\
&=&y\phi_n(xy)-\psi_n(xy),
\end{eqnarray*}

Next, replace $g(n)$ by $f(n)/((1-p)^n(1-q)^n)$, noting that $(x+1)(y+1)=1/((1-p)(1-q))$, to get the final result, namely that the recurrence for $f(n)$ is
\begin{equation}
\label{eq:frcr}
\frac{f(n+1)-f(n)}{((1-p)(1-q))^{n+1}}=y\phi_n(xy)-\psi_n(xy)=\frac{q}{1-q}\phi_n\left(\frac{pq}{(1-p)(1-q)}\right)
-\psi_n\left(\frac{pq}{(1-p)(1-q)}\right).\end{equation}

It is easy to find the generating function of the sequence $\{f(n)\}_{n=0}^{\infty}$ from the recurrence (\ref{eq:frcr}). It is
\begin{equation}\label{eq:fgf}
\sum_{n\ge 0}f(n)t^n=\frac{1}{2(1-t)}\left(1-\frac{1-(1-p+q)t}{\sqrt{(1-(1-p-q)t)^2-4pqt}}\right).
\end{equation}
\begin{corollary} (Symmetry)\label{cor:symm}
For all $p,q$ we have $f(n,p,q)=f(n,1-q,1-p)$, and consequently for all $0<q<p<1$ we have
\begin{equation}\label{eq:nsym}
N(q,p)=N(1-p,1-q).
\end{equation}
\end{corollary}
\noindent Proof 1. For a first proof, replace $p$ by $1-q$ and $q$ by $1-p$ in the generating function (\ref{eq:fgf}),
and check that it remains unchanged. $\Box$

\noindent Proof 2. For a more earthy proof, Alice's winning of the $(q,p)$ game means she had more
heads. This is identical to Bob's having more tails. That occurs when Bob wins the tails game where
he has a coin that comes up tails with probability $1-p$ and Alice has a coin that comes up tails
with probability $1-q$. The probability of the latter is $f(n, 1-p, 1-q)$ while the former happens
with probability $f(n,q,p)$. $\Box$
\subsection{Remarks on the identity (\ref{eq:fmon})}
The identity in equation (\ref{eq:fmon}) relates two different forms of the function $f(n)$, namely the form (\ref{eq:eq1}), of its original definition, and the form on the right side of (\ref{eq:fmon}). Our first comment is that although the identity was discovered by Zeilberger's algorithm, it can be given a straightforward human proof, which we will now sketch.

First, if we substitute (\ref{eq:eq1}) into the left side of (\ref{eq:fmon}) it takes the form
\begin{equation}\label{eq:left}
\sum_r\sum_{s>r}\left\{{n+1\choose r}{n+1\choose s}-(1+x)(1+y){n\choose r}{n\choose s}\right\}x^ry^s,
\end{equation}
when we express it solely in terms of the variables $x$ and $y$.

Now by inspection of the right side of (\ref{eq:eq1}) we see that only terms $x^ay^b$ appear in which $b-a=0$ or 1. Thus to prove the identity we might show that all monomials $x^ay^b$ on the left side, i.e., in (\ref{eq:left}) above, vanish if $b-a\neq 0$ or 1, and that the remaining terms agree with those on the right. We omit the details.

Our second comment is that from the identity (\ref{eq:fmon}) we can find a new formula for $f(n)$ itself, the probability that Alice wins if $n$ tosses are done. To do this, multiply both sides of (\ref{eq:fmon}) by the denominator on the left, and sum over $n$. The left side telescopes and we find
\begin{equation}\label{eq:tele}
f(n)=\frac12(1-(1-p-q)^nP_n(u)-(p-q)\sum_{k=0}^{n-1}(1-p-q)^kP_k(u)),
\end{equation}
in which, as usual, we have put $u=1+2pq/(1-p-q)$.
This formula is well adapted to computation of $f(n)$. Let's define
\[Y_n=\frac{1}{1-p-q}+\sum_{k=0}^{n-1}(1-p-q)^kP_k(u).\]
Then it's not hard to show that $\{Y_n\}$ satisfies the recurrence
\begin{eqnarray*}
(n-1)Y_n&=&(3 p+3 q-6 p q-4+n (3-2 p-2 q+4 p q))Y_{n-1}\\
&&-(7 p-2 p^2+7 q-10 p q-2 q^2-5+n(3-4 p+p^2-4 q+6 p q+q^2))Y_{n-2}\\
&&+(n-2) (1-p-q)^2 Y_{n-3},\quad (n\ge 2),
\end{eqnarray*}
with $Y_{-1}=0$, $Y_0=1/(1-p-q)$, $Y_1=1+1/(1-p-q)$. Using this recurrence in (\ref{eq:tele}) is the only way we know to compute accurate values of the probability when $n$ is large.

\section{Proof of the unimodality theorem} \label{sec:uniq}
In this section we will prove Theorem \ref{th:uni}, the unimodality theorem for the optimum value of $n$.

According to (\ref{eq:find}), we have $f(n+1)> f(n)$ precisely for those $n$ such that
$y\phi_n(xy)-\psi_n(xy)> 0$, i.e., as long as
\[\left(y+\frac12(1+xy)\right)\phi_n(xy)-\frac12 \phi_{n+1}(xy)> 0,\]
or equivalently, as long as
\[\frac{\phi_{n+1}(xy)}{\phi_n(xy)}< 1+(x+2)y,\]
or
\begin{equation}
\label{eq:ss2}
(1-xy)\frac{P_{n+1}\left(\frac{1+xy}{1-xy}\right)}{P_{n}\left(\frac{1+xy}{1-xy}\right)}< 1+(x+2)y.
\end{equation}
First suppose that $xy<1$, i.e., that $p+q<1$. We claim
\begin{theorem}
\label{th:th2}
Fix a number $x>1$.  Then the ratios
\[ \frac{P_{n+1}(x)}{P_n(x)}\qquad(n=0,1,2,\dots)\]
 strictly increase with $n$.
\end{theorem}
While Theorem \ref{th:th2} can be proved by induction on $n$, we use a more general technique which shows that the result holds not only for the Legendre polynomials, but for any sequence of functions of $n$ that are representable as $\int_a^bg(t)^nh(t)dt$, with positive $g,h$. See Lemma \ref{lem:lem1} below.

 To prove this we start with a definition and lemma.\\
 \noindent\textbf{Definition.} A function $g(t)$, defined on an interval $a\le t\le b$, is admissible for that
 interval if $g(t)\ge 0$ for all $t\in (a,b)$, and for every finite sequence $\{x_i\}$ of real numbers, not all 0,
  it is true that $\sum_ix_ig(t)^i$ does not vanish identically on $(a,b)$.

\begin{lemma}
\label{lem:lem1}
Suppose $g(t)$ is admissible for $(a,b)$, and  $h(t)\ge 0$ for all $t\in (a,b)$. Let $\mu_n=\int_{a}^bg(t)^nh(t)dt$,
 and suppose that all $\mu_i>0$.
 Then $\frac{\mu_{i+1}}{\mu_i}$ is a strictly increasing function of $i=0,1,2,\dots$.
\end{lemma}

\noindent Proof. Let $H$ be the infinite Hankel matrix $\{\mu_{i+j}\}_{i,j\ge 0}$.
Consider the principal submatrix formed by the first $n$ rows and columns of $H$.
If $x_0,x_1,\dots,x_{n-1}$ are arbitrary real numbers, not all zero, then the quadratic form
\[Q_n=\sum_{i,j=0}^{n-1}x_iH_{i,j}x_j=\sum_{i,j=0}^{n-1}x_ix_j\int_{a}^bg(t)^{i+j}h(t)dt
=\int_{a}^b\left(\sum_{i=0}^{n-1}x_ig(t)^i\right)^2h(t)dt,\]
 is clearly positive. Hence $H$ is a positive definite matrix, whence its $2\times 2$ principal minors
 $\mu_{2i}\mu_{2i+2}-\mu_{2i+1}^2$ are all positive, i.e.,
 \begin{equation}
 \label{eq:inqs}
 \frac{\mu_1}{\mu_0}<\frac{\mu_2}{\mu_1};\ \frac{\mu_3}{\mu_2}<\frac{\mu_4}{\mu_3};\dots
\end{equation}
Next replace $h(t)$ by $g(t)h(t)$. Then the sequence $\{\mu_i\}_{i\ge 0}$ is replaced by $\mu_{i+1}$ (${i\ge 0}$) and
the Hankel matrix $H$ is replaced by one whose $(i,j)$ entry is $\mu_{i+j+1}$. We apply the conclusion (\ref{eq:inqs})
 to this new situation and we discover that
 \begin{equation}
 \label{eq:inqs2}
 \frac{\mu_2}{\mu_1}<\frac{\mu_3}{\mu_2};\ \frac{\mu_4}{\mu_3}<\frac{\mu_5}{\mu_4};\dots
\end{equation}
If we combine (\ref{eq:inqs}) and (\ref{eq:inqs2}) we obtain
the result stated in Lemma \ref{lem:lem1}. $\Box$
\bgs
To prove Theorem \ref{th:th2} we have the integral representation
\begin{equation}\label{eq:intrep}
P_n(x)=\frac{1}{\pi}\int_0^{\pi}(x+\sqrt{x^2-1}\cos{t})^ndt
\end{equation}
for the Legendre polynomials. We can take $g(t)=x+\sqrt{x^2-1}\cos{t}$ and $h(t)=1$ in Lemma \ref{lem:lem1} and the
conclusion of Theorem \ref{th:th2} follows. $\Box$
\bgs
 We remark
 that this is the reversal of a
  celebrated inequality of Tur\' an which holds inside the interval of orthogonality. Now that the ratio of
  the Legendre polynomials on the left side of (\ref{eq:ss2}) is known to be a strictly increasing function
   of $n$, we observe that when $n=0$ the left side has the value $1+xy$, and when $n\to\infty$, the well
    known asymptotic behavior of $P_n(x)$ for fixed $x>1$ and large $n$ shows that the left side
    approaches $(1+\sqrt{xy})^2$, which is larger than $1+(x+2)y$. Hence there is a unique $n$ for which the
left side of (\ref{eq:ss2}) is $\le$ the right side, but at $n+1$ the inequality is reversed.

The case where $xy>1$ is reduced to the case $xy<1$, which we have just handled, by equation (\ref{eq:nsym}). If $xy=1$, i.e., if $p+q=1$,
we discuss the situation in the next section. The proof of Theorem \ref{th:uni} is now complete. $\Box$

\section{The interesting special case $p+q=1$}

Consider the special case in which $p+q=1$. Then $xy=p(1-p)/((1-p)p)=1$, and we can carry out the
calculations analytically in full. Indeed, we now have
$\phi_n(1)={2n\choose n}$ and $\psi_n(1)=n{2n\choose n}/(n+1)$, from which we get
\begin{eqnarray}\label{eq:jqq}
J_n(q,1-q)&=&\frac{q}{1-q}\phi_n(1)-\psi_n(1)=\frac{q}{1-q}{2n\choose n}-\frac{n}{n+1}{2n\choose n}\nonumber\\
&=&\left(\frac{q}{1-q}-\frac{n}{n+1}\right){2n\choose n}.
\end{eqnarray}
 This last vanishes iff $q/(1 - q )= n/(n+1)$, which is $q=n/(2n+1)$. The sign of $J_n$ then equals that of $q-n/(2n+1)$. This proves the following.
\begin{lemma}\label{lem:diagj}
\[J_n\left(\frac{n}{2n+1},1-\frac{n}{2n+1}\right)=0.\]\end{lemma}
Now the unimodality theorem, Theorem \ref{th:uni}, gives an explicit formula for $N$ on this line.
\begin{theorem}(The Diagonal Formula)\label{th:diagf}
   For $0<q<1/2$, we have  \[ N(q,1-q)=\left\lceil{\frac{q}{1-2q}}\right\rceil.\]
\end{theorem}
\noindent Proof.  By the uniqueness theorem, and the tie-breaking aspect of the definition of $N(q,p)$, $N(q,p)$ is always the least integer $n$ such that $J_n\le 0$. Because of the agreement of the sign of $J_n$ and that of $q-n/(2n+1)$, the result follows. Note that when $q=n/(2n+1)$, the $J=0$ condition means that there is a tie between
the two values $\lceil{q/(1-2q)}\rceil$ and $\lceil{q/(1-2q)}\rceil+1$ for the optimal choice. $\Box$

\section{A general lower bound}
\begin{theorem}\label{th:low} Let $N(q,p)$ denote the optimum choice of $n$, i.e., the one that maximizes Alice's chance of winning. Then
we have
\begin{equation}
N(q,p)\ge \left\lfloor{\frac{1}{2(p-q)}-\frac12}\right\rfloor .
\end{equation}
\end{theorem}
First we need the following
\begin{lemma}\label{lem:trap}
 In the trapezoid $\tau_n$, defined by the lines $q=0$, $p=q$, $p+q=1$, and the inequality $p\le 1/(2n+1)+q$, the indicator function $J_n$ is positive. That is, if $p\le L_n(q)$ and $(q,p)\in T$, then $J_n(q,p)>0$.
\end{lemma}
\noindent Proof of the Lemma. We use induction on $n$. Figure 2 shows the trapezoid.  Since $\tau_n\subseteq \tau_{n-1}$, the induction is valid.
\begin{figure}[htp]\label{fig:fig2}
\begin{center}
\includegraphics{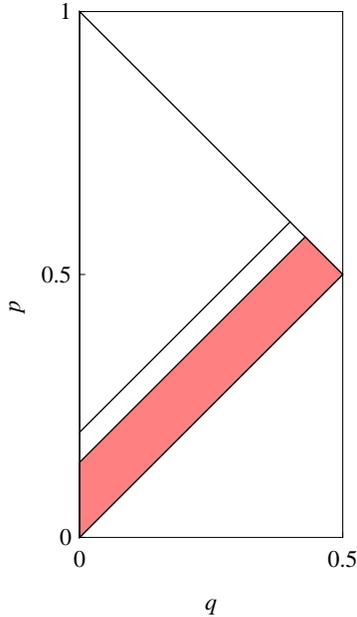}
\caption{The trapezoid $\tau_n$}
\end{center}
\end{figure}
The case $n=1$ follows easily from
$J_1(q,p)(1-p-q)/(2 q (1-p))=1-q+p (-2+3 q),$
$L_1=1/3+q$, and $q<1/2$. Suppose that $J_n\le 0$. Because $p<L_{n-1}(q)$, the induction hypothesis applies, giving $J_{n-1}>0$. Therefore
$r_{n+1}\ge \rho>r_n$. The ratio recurrence tells us that $(n+1)r_{n+1}-(2 n+1) u+n/ r_n=0.$ Therefore $(n+1)\rho-(2 n+1) u+n/\rho<0$,
or $1-(2 n+1) p+2 n q+q<0$. Thus $1/(2 n+1)+q<p$, contradicting $p\le L_{n}(q)$. $\Box$
\bgs
\noindent\textbf{Proof of Theorem \ref{th:low}.} Again, suppose first that $xy<1$, i.e., that $p+q<1$. The theorem says that if $p\le L_n(q)$, then  $n$ cannot be $N(q,p)$. Therefore for $n\le 1/(2(p - q)) - 1/2$, $n$ is not $N(q, p)$ and $N(q, p) > 1/(2(p - q))- 1/2$. More precisely,
$N\ge \lceil{1/(2 (p - q)) - 1/2}\rceil = \lfloor{ 1/(2 (p - q)) + 1/2}\rfloor$.     But on the $p+q=1$ line, \[N(q,p)=\left\lceil{q/(1-2q)}\right\rceil=\lfloor{1/(2(p-q))-1/2}\rfloor.\]  So the latter works as a lower bound in both cases. Finally, if  $p+q>1$, then the symmetry formula of Corollary \ref{cor:symm} yields
$N(q,p) = N(1-p,1-q)$, a transformation that leaves the bound invariant. $\Box$

\section{The upper bound}\label{sec:upb}
In this section we study in detail the curves $J=0$ and use the results to obtain a simple upper bound on $N(q,p)$ which is  roughly twice
 the lower bound of Theorem \ref{th:low}.

\begin{theorem}(Upper bound on $N$)\label{th:ubn}
\[N(q,p)\le \frac{\max{(1-p, q)}}{p-q}.\]
\end{theorem}
\subsection{Curves on which $J_n$ vanishes}
The key idea underlying our analysis of $N(q,p)$ is an understanding of the vanishing sets of $J_n$. Figure 3 shows these curves, together
with some of the lines $M_n$ and $L_n$.
\begin{figure}[htp]\label{fig:fig3}
\begin{center}
\includegraphics{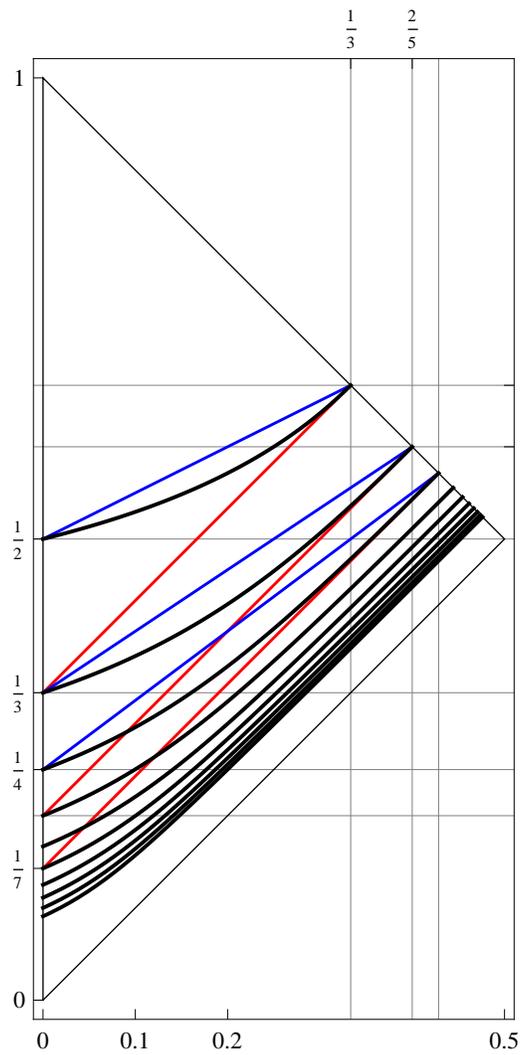}
\caption{The black curves are where $J_n(q,p)$ vanishes. The blue lines are $M_n$ and the red lines are $L_n$, for $n=1,2,3$. }
\end{center}
\end{figure}

The upper black curve is $J_1=0$, and so on down. The blue lines are $M_n$ and the red ones, $L_n$. So we know that below the uppermost
 black curve $J_1>0$ and so $N\ge 2$. In fact, the $N=1$ region is just the region above the first black curve. But now we need to prove
 various properties evident from the diagram.
\subsection{The function $p_n(q)$}
Our first task will be to define, and to verify the correctness of the definition, of a function $p_n(q)$ which for each $q\in (0,n/(2n+1)]$,
 is the unique value of $p$ for which $J_n(q,p)=0$. We will therefore start by proving existence and uniqueness of such a $p$.
\begin{lemma}\label{lem:lem2}
For $q\in T_n$, $J_n(q,1-q)<0$.
\end{lemma}
\noindent\textbf{Proof.} As in eq. (\ref{eq:jqq}) we have
\[J_n(q,1-q)=\left(\frac{q}{1-q}-\frac{n}{n+1}\right){2n\choose n},\]
which is negative iff $q<n/(2n+1)$. $\Box$
\begin{lemma}\label{lem:lem3}
For $q\in T_n$, $J_n(q,q)>0$.
\end{lemma}
\noindent\textbf{Proof.} The condition that $J_n(q,q)>0$ is the same as $r_{n+1}<\rho=1/(1-2q)$, so we must prove that $r_{n+1}<1/(1-2q)$ when $q\in T_n$. We will prove more, namely that $r_n<1/(1-2q)$ whenever $0<q<1/2$.

Let $r$ be the fixed point that is $>1$ of the Legendre polynomial recurrence, i.e. the root of the quadratic equation
\[r=\frac{(2 n+1)u}{n+1}-\frac{n}{(n+1)r}\]
that is $>1$. This root is
\[W_n=\frac{\sqrt{(2 n+1)^2 u^2-4 n(n+1)}+u(2 n+1)}{2 (n+1)}.\]
It is easy to check that $1<W_n\le 1/(1-2q)$, and tedious, but routine, to check that $W_n<W_{n+1}$ for $0<q<1/2$.

Now we can prove by induction that $r_n<W_n$ whenever $0<q<1/2$, which suffices. (Note the change from $q\in T_n$ to $q\in (0,1/2)$; this is essential to allow the induction to carry through.) The base case can be taken to be
\[W_1-r_1=\frac{1}{8} \left(2 q+2
   \sqrt{-9 (1-q) q+\frac{9}{4
   (1-2
   q)^2}-\frac{5}{4}}-\frac{1}
   {1-2 q}-1\right)\]
whose positivity is follows from the fact that the expression is positively infinite at $q=1/2$ and $0$ only when $q=0$.

For the inductive step, take the recurrence
\[r_{n+1}=\frac{(2n+1)u}{n+1}-\frac{n}{(n+1)r_n}\]
 and assume inductively that $r_n<W_n$. Then
 \[r_{n+1}<\frac{(2 n+1)u}{n+1}-\frac{n}{(n+1) W_n},\] but this last equals $W_n$, a fixed point of the recurrence. Therefore $r_{n+1}<W_n$,
 and the proof of Lemma \ref{lem:lem3} is complete because $W_n<W_{n+1}$ for $0<q<1/2$. $\Box$

 The algebraic part of the proof actually yields the more general result that $r_n < 1/(1-p-q)$.

\begin{theorem}\label{th:exist} (Existence theorem) Given $n\ge 1$ and $q\le n/(2n+1)$, there is a value $p$ with $q<p\le 1-q$ such that  $J_n(q,p)$ vanishes.
\end{theorem}
\noindent{Proof}. When $q<n/(2(n+1))$ this follows from Lemmas \ref{lem:lem2} and \ref{lem:lem3}, which tell us that $J_n(q,p)$
is a strictly decreasing function of $p$, going from a positive value to a negative one. At the endpoint the fact that $J_n(q,1-n/(2(2n+1))=0$ gives the existence of the desired $p$. $\Box$
\bgs
\noindent \textbf{Definition}: Given $n\ge 1$ and $q\in T_n$, we write $p_n$ for the largest value of $p$ between $q$ and $1-q$ such that $J_n(q,p)=0$.

We now turn to the important proof that, in all cases, there is only one value $p$ so that $J_n(q, p) = 0$. The key is the following lemma about $\partial J/\partial p$.
\begin{lemma}\label{lem:lem4} We have
\[\frac{\partial}{\partial p}J_n(q,p)<0\  \mathrm{for}\ q\le p\le 1-q.\]
\end{lemma}
\noindent Proof. The derivative inequality, after multiplication by
\[ \frac{2}{q} (1-p)^3 p (1-q)^2\left(\frac{-p-q+1}{(1-p)(1-q)}\right)^{1-n},\]
becomes
\[P_n(u) \left(n \left(2 p^2-p(2 q+1)+q-1\right)-2 pq+p+q-1\right)+(n+1)(-p-q+1) P_{n+1}(u)<0,\]
and so holds precisely when $r_{n+1}<V$, where
\[V=\frac{n \left(2 p^2-p (2
   q+1)+q-1\right)-2 p
   q+p+q-1}{(n+1) (p+q-1)}.\]
(The
inequality was reversed because $-(n+1) (1-p-q)<0$.)

Now $r_{n+1}<V$ can be proved by an easy induction, since the domain of truth does not depend on $n$. For the base case examine $V-r_1$, which works out to be the positive quantity
\[\frac{2 np (1-p) }{(n+1) (1-p-q)}.\]
 The induction step uses the usual recurrence. Suppose $r_n<V$; then
 \[r_{n+1}=\frac{(2 n+1) u}{n+1}-\frac{n}{(n+1) r_n}<\frac{(2
   n+1) u}{n+1}-\frac{n}{(n+1)V}.\]
    But this last is less than $V$ because the difference
    $V-((2 n+1)u)/(n+1)-n/((n+1) V)$ works out to
\[\frac{2 n (p-q)-2 q+1}{n \left(-2 p^2+2 p
   q+p-q+1\right)+2 p q-p-q+1}\] in which all the grouped terms are positive. $\Box$

\begin{theorem}(Uniqueness Theorem)\label{th:uniq} Given $n\ge 1$ and $0<q\le n/(2n+1)$. If $q\le p<p_n(q)$, then $J_n(q,p)>0$. Thus there is only one vanishing value for $J_n$.
\end{theorem}
\noindent Proof.  This follows from Lemma \ref{lem:lem4} which implies that $J_n(q,p)$ is a strictly
 decreasing function of $q$; thus it cannot return to 0 after it has once taken that value. $\Box$

\subsection{Properties of $J$}
 Now we have the functions $p_n$ defined on $T_n$, with the property that $J_n(q,p_n(q))=0$. We will call the graph of $p_n(q)$ a $J$-\textit{nullcline}, and denote it simply by $p_n$, or often just $p$. We need several properties of $p_n$. Note that most of the properties below have one-line proofs thanks to the efficient definition of $p_n$ and the uniqueness result (Theorem \ref{th:uniq}).

\begin{lemma}\label{lem:fact2} The function $p_n$ is continuously differentiable ($C^1$) on $T_n$.\end{lemma}
\noindent Proof. This is a consequence of the fact that $J_n(q,p)$ is differentiable (it is a rational function) and Lemma \ref{lem:lem4}, which states that $\frac{\partial}{\partial p}J_n>0$. These facts show that the hypotheses of the implicit function theorem are satisfied.

The preceding result is about the triangle $T$. But it also works on the $p+q=1$  line, for on that line
\[J_n(q,p)=J_n(q,1-q)= \left(\frac{q}{1 - q }- \frac{n}{n+1}\right){2n\choose n}=\left(\frac{q}{p}-\frac{n}{n+1}\right){2n\choose n},\]
 and the partial derivative $\frac{\partial}{\partial p}J_n(q,p)$ is just $-q{2n\choose n}/p^2$,
 which is nonzero. By symmetry the same proof works on the opposite side of the $p+q=1$ line. $\Box$

\begin{lemma}(Derivative formula)\label{lem:fact3}
For any point $(q,p)\in T$ and on the graph of $p_n$, we have
\[ p_n'(q)=\frac{p(1-p)(np-nq+p-1)}{q(1-q)(n(q-p)+q)}.\]
\end{lemma}
\noindent Proof. By the implicit function theorem,
\begin{equation}\label{eq:ift}
p_n'=-\frac{\frac{\partial}{\partial q}J_n(q,p)}{\frac{\partial}{\partial p}J_n(q,p)}.
\end{equation}
Taking the derivatives, using the recurrence to eliminate $P_{n+2}$, and using the
$J_n(q,p)=0$ relation to eliminate $P_{n+1}$ in favor of $P_n$ leads immediately to the formula. $\Box$
\begin{lemma}(Linear vanishing condition)\label{lem:fact4}
 $p_n'(q)=0$ iff $p_n$ lies on the line $M_n$.
 \end{lemma}
\noindent Proof: Immediate from the numerator of the derivative formula (Lemma \ref{lem:fact3}). $\Box$

\begin{lemma}\label{lem:fact1}
Given $n\ge 1$ and $q<p$, if $J_n(q,p)=0$ then the partial derivative $(\partial/\partial q)J_n(q,p)$ is not zero.\end{lemma}
\noindent Proof.
We work first in $T$. The partial derivative of $J_n(q,p)$ can be taken and simplified using the standard recurrence for $P_{n+2}$, and also the relationship derived from $J_n=0$, to replace $P_{n+1}$ by $\rho P_n$. This leads to
\[\frac{\partial}{\partial q}J_n(q,p)=-\frac{P_n(u) (n (p-q)+p-1)
   \left(\frac{-p-q+1}{(p-1)
   (q-1)}\right)^n}{(1-q)^2
   (-p-q+1)}.\]
The denominator is nonzero and $xy\ge 1$ in $T$ so $P_n(u)\ge 1$; therefore the partial derivative vanishes at a point on a $J$-nullcline iff $n(p-q)+p-1=0$ iff
$p=\frac{1}{n+1}+\frac{n}{n+1} q=M_n(q).$
 So we must show that this value of $p$ cannot lead to a point at which $J_n$ vanishes. Define $h(q)=J_n(q,M_n(q))$. This evaluates to
\[\frac{\left(\frac{-2 nq+n-q}{n (q-1)^2}\right)^n
   ((n+q) P_n(u)+(n (2 q-1)+q)P_{n+1}(u))}{2 n (1-q)^2}.\]
When $q=0$, we have $u=1$ and so $P_n(u)=1$ and $h (0)=0$. We also have
(limits are from the left)
\[\lim_{q\to n/(2n+1)}h(q)=\lim_{q\to n/(2n+1)}J_n(q,M_n(q))=J_n\left(\frac{n}{2n+1},\frac{n+1}{2n+1}\right),\]
by continuity of $J_n$. But this last vanishes, by the diagonal formula of Theorem \ref{th:diagf}. Thus
\[\lim_{q\to n/(2n+1)}h(q)=0.\] Hence $h(q)$ is a differentiable function of $q$ which vanishes at both ends of the interval $(0,n/(2n+1))$. Figure 4 shows the graphs of $h$ for $n=1,2,3,4,5$. It remains to show that $h$ cannot vanish for any $q\in T_n$.

\begin{figure}[htp]\label{fig:fig4}
\begin{center}\includegraphics{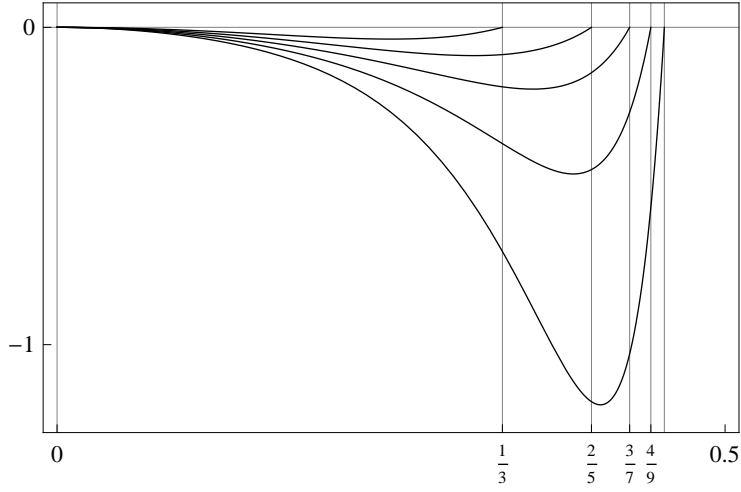}\end{center}
\caption{The graphs of $J_n(q,M_n(q))$ for $n=1,2,3,4,5$.}
\end{figure}

Suppose $h(q)$ vanishes at some $q\in T_n$. Then, because of the vanishing at the endpoints, there must be a point $q_1$ such that $h'(q_1)=0$ and $h(q_1)\ge 0$. More precisely, if there is a true crossing at $q$ then $q_1$ would be given by the Mean Value Theorem applied to either $[0,q]$ or $[q,n/(2n+1)]$; if there is a tangency to the axis at $q$ (or even if the function is identically 0), then $q_1=q$.

So at the point $(q_1,M_n(q_1))$ (which we denote by just $(q,M_n(q))$ in the expressions below) we would have two relations for $r_{n+1}$: The $h(q_1)=J_n(q_1, M_n(q_1))\ge 0$ condition means that
\[ r_{n+1}\le \rho=(1-p+q_1)/(1-p-q_1),\]
 where in the last fraction $p$ is to be replaced by $M_n(q_1)$; and taking the $q$-derivative of $h(q)=J_n(q,M_n(q))$ and recursively removing $P_{n+2}$ leads to the following equation:
\[-(n+1)\left(\frac{n-q-2nq}{n(q-1)^2}\right)^n\frac{c_1P_n(u)+c_2P_{n+1}(u)}
{2n(q-1)^3(nq+1)(n(2q-1)+q)}=0,\]
where
\[c_1=2 n^3 q+n^2 \left(4 q^3-4 q^2+5 q+1\right)+n
   \left(2 q^3+q^2+2 q+1\right)+q (q+1),\]
\[c_2=(n+1)  (1+q+2 n q) (-n+q+2 n q).\]
Clearing the nonzero factors tells us that $r_{n+1}=-c_1/c_2$, and therefore
\[-\frac{c_1}{c_2}\leq \frac{-\frac{n
   q}{n+1}-\frac{1}{n+1}+q+1}{-\frac{n
   q}{n+1}-\frac{1}{n+1}-q+1}\]
   which reduces to
   \[\frac{2 n (1-q) q}{(n+1) (2 n q+q+1)}\leq 0,\]
   a clear contradiction, which establishes the theorem for the triangle $T$. But the same proof works on the $p+q=1$ line, for on that line \[J_n(q,p)=J_n(q,1-q)= \left(\frac{q}{1-q}-\frac{n}{n+1}\right){2n\choose n}=\left(\frac{q}{p}-\frac{n}{n+1}\right){2n\choose n},\]
and the partial derivative $(\partial /\partial q)J_n(q,p)$ is just ${2n\choose n}/p$
 which is nonzero. By symmetry the same proof works on the opposite side of the $p+q=1$ line. $\Box$
\begin{lemma}\label{lem:fact5}
For $q\in T_n$ we have $p_n'(q)\neq 0$.
\end{lemma}
\noindent Proof: Because of (\ref{eq:ift}),
 the claim follows from Lemma \ref{lem:fact1},
  which shows that the numerator does not vanish on the graph of $p_n$. $\Box$
\begin{lemma}\label{lem:fact6}
\[p_n\left(\frac{n}{2n+1}\right)=1-\frac{n}{2n+1}.\]
\end{lemma}
\noindent Proof: Follows from the diagonal equation of Lemma \ref{lem:diagj}. $\Box$
\begin{lemma}\label{lem:fact7}
\[p_n'\left(\frac{n}{2n+1}\right)=1.\]
\end{lemma}
\noindent Proof: The definition of $p_n$ can be carried over by symmetry to the other side of the $p+q=1$ line to yield a differentiable function. The implicit function theorem applies on the line itself as noted in the proof. But then, by symmetry of all the probabilities, and therefore of the vanishing of $J_n$, $p_n$ is symmetric
across the line. Thus differentiability implies that the derivative must be 1 to avoid a cusp. $\Box$
\begin{theorem}(Upper bound on $p_n$)\label{th:upper} For every $q\in T_n$, we have $p_n(q)<M_n(q)$, where $M_n(q)$ is the line (\ref{eq:mnq}) above. \end{theorem}
\noindent Proof: Because $p_n$ and $M_n$ agree on the line $p+q=1$ (Lemma \ref{lem:fact6}), and because \[p_n'\left(\frac{n}{2n+1}\right)=1,\]
 (Lemma \ref{lem:fact7}) while the slope of $M_n$ is $n/(n+1)$, the fact that $p_n$  is $C^1$ (Lemma \ref{lem:fact2}) means that $p_n$ is under $M_n$ when $q$ is just left of the $p+q=1$ line. Lemma \ref{lem:fact5} tells us that $p_n'$ is never 0, and therefore, by Lemma \ref{lem:fact4},  $p_n$ can never cross the line $M_n$. $\Box$
\begin{lemma}\label{lem:fact9}
 $N(q,p_n(q))=n$.
\end{lemma}
\noindent Proof.  When $(q,p)\in T$, this follows from Theorem \ref{th:uni} because $J_n(q,p_n(q))=0$ and so $J_{n-1}(q,p_n(q))>0$. Therefore $n$ is the least $m$ such that $J_m(q,p_n(q))=0$. If $p_n(q)=1-q$ then it must be that $q=n/(2n+1)$ and, by the proof of Theorem \ref{th:uni},
 $n=q/(1-2q)$ is the unique integer such that $J_n(q,1-q)=0$, establishing the result. $\Box$
\begin{lemma} \label{lem:fact12} For every $q$ we have $p_1(q)>p_2(q)>p_3(q)>\dots$.\end{lemma}
\noindent Proof. The graphs $p_n$ can never cross because if $p_n(q)=p_m(q)$, $N (q,p)$ would be simultaneously $n$ and $m$ by Lemma \ref{lem:fact9}. and the right end of $p_n$ is above the
right end of $p_{n+1}$ (Lemma \ref{lem:fact6}). So continuity of the graphs yields the result. $\Box$

\begin{lemma}\label{lem:lem13}
The graphs of $p_n(q)$ determine the value of $N(q,p)$ exactly as follows: $N(q,p)$ is the least $n$ such that $p\ge p_n(q)$.
\end{lemma}
\noindent Proof. We know this is correct when we are on the graph $p_n$ (Lemma \ref{lem:fact9}). But if $(q,p)$ is between $p_{n-1}$
 and $p_n$ then $J_{n-1}>0$ and $J_n<0$, by Theorem \ref{th:exist}, and this means $N(q,p)=n$. $\Box$
\bgs
Figure 5 shows how the graphs of $p_n$ divide triangle $T$ into regions that define the optimal $N$-values, and how $p_n$ is bounded by the two lines $M_n$ and $L_n$.

\begin{figure}[htp]\label{fig:fig5}
\begin{center}
\includegraphics{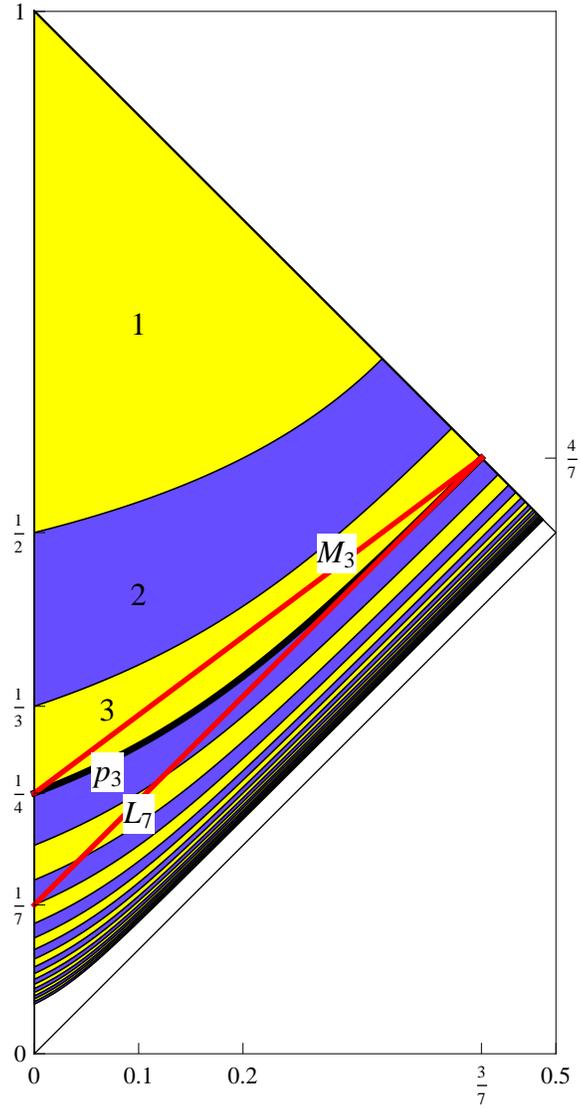}
\caption{The graphs of $p_n$ separate $T$ into regions where $N=1,2,3,\ldots.$ The two red lines $M_3$ and $L_3$ form bounds on $p_3$, and this relationship underlies the lower and upper bounds on $N(q,p)$.}
\end{center}
\end{figure}

\begin{corollary}\label{cor:big}
$N(q,p)$ is a nonincreasing function of $p$. Thus if the $p$-coin becomes stronger,
 then the optimal choice of game length for the holder of the $q$-coin cannot get larger.
\end{corollary}

\noindent We can now prove Theorem \ref{th:ubn}, the upper bound on $N(q,p)$.\\
\noindent Proof. Assume first that $(q,p)\in T$, in which case the claimed bound is just $\left\lceil{(1-p)/(p-q)}\right\rceil$. Let $n$ be the smallest value so that $p$ lies at or above $M_n(q)$. Then $n=\left\lceil{(1-p)/(p-q)}\right\rceil$ and for this to be a bound we need, by Lemma \ref{lem:lem13}, that $p_n(q)<M_n(q)$. But this is exactly what Theorem \ref{th:upper} tells us. For the special case on the diagonal line: $N(q,1-q)=\left\lceil{q/(1-2q)}\right\rceil$, by the diagonal formula. But this is identical to $\left\lceil{(1-p)/(p-q)}\right\rceil$, which therefore works in both cases. The case in which $p+q>1$
 is handled by symmetry (see Corollary \ref{cor:symm}), with $q$ taking the role of $1-p$. $\Box$
\bgs
The two proved bounds in Theorem \ref{th:bounds} are equal 47\% of the time because they must agree if it happens that $M_n(q) < p < L_{n-1}(q)$; these conditions define a collection of triangles whose area, in proportion to $T$, is easily computed to be $\pi^2/4-2$. In all such cases, then, $N (q,p)$ equals $\lceil{(1-p)/(p-q)}\rceil$.
\section{Deeper analysis of the nullclines}\label{sec:sec9}
In section \ref{sec:upb} we proved many properties of $p_n$ that were evident from the graphs. We continue that here, gaining information that leads to improved bounds on $N$ and to an efficient algorithm for computing $N$ (section \ref{sec:sec10}). We first observe that when $n$ is small, $p_n$ is given by simple formulas. Such formulas are useful as a check on computations.
\begin{lemma}\label{lem:small}
\begin{eqnarray*}
p_1(q)&=&\frac{1-q}{2-3q},\\
p_2(q)&=&\frac{1-q}{3-12q+10q^2}(2-4q-\sqrt{1-4q+6q^2}).
\end{eqnarray*}
\end{lemma}
\noindent Proof. Just solve the polynomial equation $J_n(q,p)=0$. There are more complicated radical expressions for $p_2$, $p_3$, and $p_4$. The last is ostensibly a
 quintic, but the polynomial in the equation is divisible by $1-p$ yielding a quartic equation.  $\Box$
 \bgs
Because $J_n$ is infinitely differentiable in $T$, the fact that the hypothesis of the implicit function theorem is met (Lemmas \ref{lem:lem4} and \ref{lem:fact2}) means that $p_n$ is infinitely differentiable. The second derivative is easily computed by implicitly differentiating the derivative formula (Lemma \ref{lem:fact3}).
\begin{lemma}\label{lem:xx1} We have \[p_n''(q)=\frac{p(1-p)(1-p-q)}{q^2(q-1)^2(n(p-q)-q)^3}Z,\]
where $Z$ is given by
\begin{eqnarray*}
&&2n^3(p-q)^3+2n^2(2 p^3-p^2 (7 q+1)+p q (7 q+3)
-2 q^2 (q+1))+n (2 p^3-p^2 (11 q+2)\\
&&\qquad\qquad +p q (11 q+10)-q (2 q^2+7 q+1))-(p-1) q (3 p-3 q-1).
\end{eqnarray*}
\end{lemma}

But we can prove the weaker and still very useful assertion that the slope never rises beyond its value at the right end.
\begin{lemma}\label{lem:xx2}
For $q\in T_n$ we have $p_n'(q)>0$.
\end{lemma}
\noindent Proof. Because the slope is 1 at the right end (Lemma \ref{lem:fact7}), but never vanishes (Lemma \ref{lem:fact5}), it is always  positive. $\Box$
\bgs

We next move to a proof of the (computationally evident) fact that the $J$-nullclines are convex. Figure 6 shows the second derivatives of $p_n$ for $n\le 15$. They are evidently unbounded at $q=0$, converging to 0 at the right, and always positive. We will now prove positivity (i.e., convexity of the $J$-nullcline curves), which will be important as a source of new approximations to $N (q,p)$.
\begin{figure}[htp]\label{fig:fig6}
\begin{center}
\includegraphics{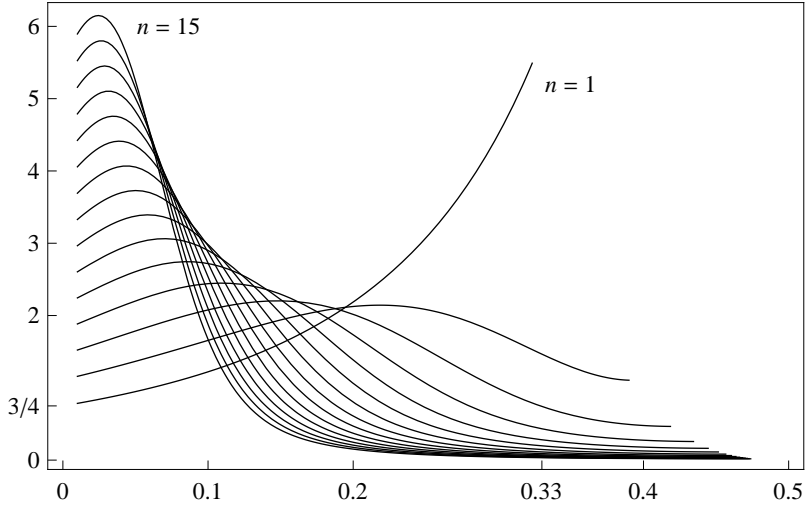}
\end{center}
\caption{The graphs of $p_n''(q)$ for $n\le 15$.}
\end{figure}
\begin{theorem}\label{th:cvx} (Convexity) $p_n''(q)>0$ in $T_n$.\end{theorem}
\noindent Proof. By Lemma \ref{lem:trap}, any point $(q,p_n(q))$ lies above $L_n(q)$, so the line connecting $(q,p_n(q))$ to $(n/(2n+1),1-n/(2n+1))$ has slope less than 1. By the Mean Value Theorem, there is some $q_1>q$ for which $p_n'(q)<1$. But the slope is 1 at $n/(2n+1)$, so there is a point at which the second derivatgive is positive. Since $p_n''$ is continuous, the proof will be complete once we show that the second derivative cannot vanish. The second derivative is given by the formula of Lemma \ref{lem:xx1}. We can eliminate nonzero factors, thus reducing its vanishing to the following equation in $p$.

\begin{eqnarray*}
&&\left(2 n^3+4 n^2+2n\right) p^3+p^2 \left(-6n^3 q-2 n^2 (7q+1)-n (11q+2)-3 q\right)\\
&&+p \left(6n^3 q^2+2 n^2 q (7 q+3)+n q(11 q+10)+3 q^2+4q\right)-n q \left(2 q^2+7q+1\right)\\
&&\qquad -3 q^2-q -2 n^3 q^3-4 n^2 q^2(q+1)=0
\end{eqnarray*}
Fixing $n$ and $q$, the vanishing condition is a cubic in $p$. The cubic always has a real root and checking the critical points assuming $0<q<1/2$ and $n\ge 1$ shows that they never straddle 0, which means that the other two roots are never real. Call the unique real root $p_n^-(q)$. The theorem will be proved once we show that $J_n(q,p_n^-(q))\neq 0$ in $T_n$, so that the inflection point is not on the graph of $p_n$. We can get a closed form for $p_n^-(q)$ by solving the cubic, but we do not need the explicit representation as the needed algebra can be worked out implicitly from the cubic relation. Yet it is instructive to look at $p_n^-(q)$. Figure 7 shows its graph for $n\le 5$, together with the graphs of $p_n$ in pink, and also the base-10 log of the difference for $n=1,2,3,10,20,100$. The inflection curve is barely below the $J$-nullcline and the two curves are visually indistinguishable.
\begin{figure}[htp]\label{fig:fig7}
\begin{center}
\includegraphics[width=6.5in]{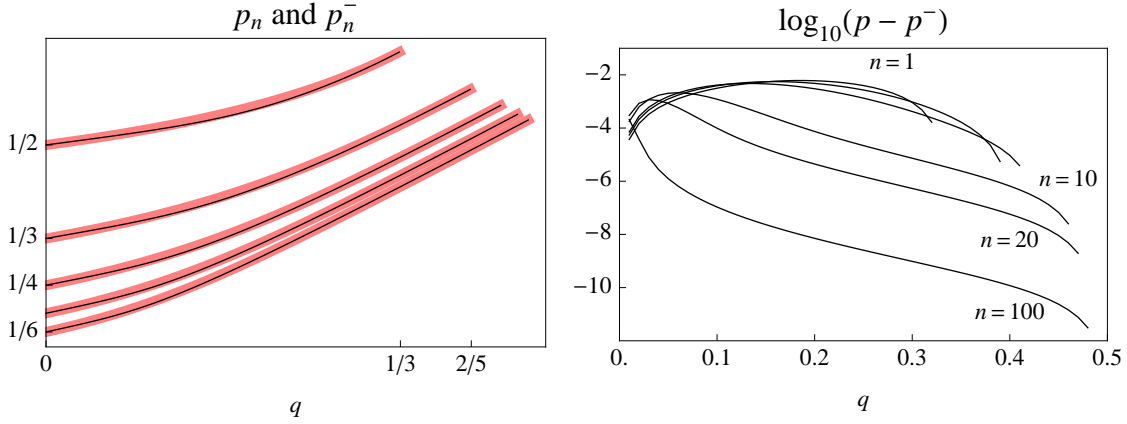}
\end{center}
\caption{Left: The graphs of $p_n$ (pink) and $p_n^-(q)$ (black). Right: The difference in the graphs viewed through a logarithmic lens.}
\end{figure}

We can use implicit differentiation on the defining cubic to get a formula for $\frac{d}{dq}p_n^-(q)$. Since $p_n^-(q)$ is given explicitly by radicals we know its derivative exists. The implicit derivative formula $p_n^-(q)=-\partial_q/\partial_p$ then gives the following representation for the derivative  $\frac{d}{dq}p_n^-(q)$:
\begin{equation}\label{eq:derv}
\frac{\left(6 n^2+8 n+3\right)p^2-2 p \left(6 n^2 q+n (8q+3)+3 q+2\right)+6 n (n+1)q^2+(8 n+6) q+1}{6 n^2(p-q)^2+2 n \left(3 p^2-2 p(4 q+1)+q (4 q+3)\right)+q(-6 p+3 q+4)}.
\end{equation}

Now we follow the proof idea of Lemma \ref{lem:fact1}. We need to show that, given $n$, the point $(q,p_n^-(q))$ cannot lie on the $J_n$-nullcline. Define the function $h (q)=J_n(q,p_n^-(q))$, which we claim does not vanish in $T_n$. Figure 8 shows the first few graphs of $h$.
\begin{figure}[htp]\label{fig:fig8}
\begin{center}
\includegraphics{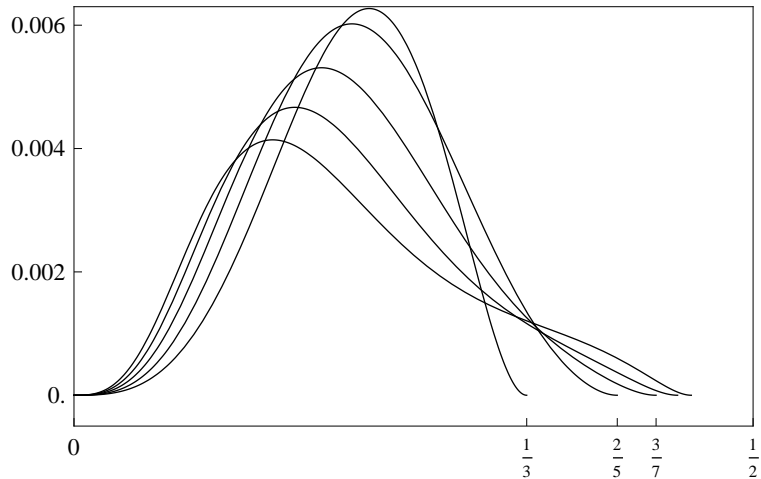}
\end{center}
\caption{The graphs of $J_n(q,p_n^-(q))$ for $n\le 5$.}
\end{figure}

We observe that $h$ vanishes at both ends. At the left this is because $J_n(0,p)$ is always 0 (the Legendre terms become just 1). At the right we find that the defining cubic, when $q$ is set to $n/(2n+1)$, has the factor $2 n p-n+p-1$, which means that $p=1-n/(2n+1)$ is a root. The cubic has exactly one real root in the domain of interest, so $1-n/(2n+1)=p_n^-(\frac{n}{2n+1}))$ and we know (Lemma \ref{lem:diagj}) that $J_n((n/(2n+1),1-n/(2n+1))=0$.

Now suppose the graph is $0$ at some $q_1\in T_n$. Then looking left and right of $q_1$ and using the Mean Value Theorem, we get a value, call it $q$, such that $h (q)\le 0$ and $h'(q)=0$. Recall that $h(q)\le 0$ means that $r_{n+1}=P_{n+1}/P_n\ge \rho=(1-p+q)/(1-p-q)$, where $p$ denotes $p_n^-(q)$. When we form $h'(q)$ leaving $p_n^-(q)$ undefined and then substitute the derivative formula (\ref{eq:derv}) for the derivatives of $p_n^-(q)$  that appear, we obtain the familiar form $c_1P_n+c_2P_{n+1}=0$. This becomes $-c_1/c_2=r_{n+1}$, so we have the relation $-c_1/c_2\ge \rho$. If we work out $c_1,c_2$ in terms of $n,p,q$ we get a rational function of these three variables which is too large to reproduce here, but which can be found at \cite{ws}. A call to Mathematica's \texttt{Reduce} function shows, in two seconds, that the five conditions:
\begin{enumerate}
\item $-c_1/c_2 \ge \rho$,
\item the vanishing of the second derivative formula at $p$,
\item $0<q<p<1-q$,
\item $n\ge 1$, and
\item $q<n/(2(n+1))$
\end{enumerate}
 are contradictory. For such polynomial systems \texttt{Reduce} uses a cylindrical algebra decomposition \cite{as}; this example requires showing that there is no solution in each of 1062
  cylindrical cells. Working this out by hand might be extremely difficult, if not impossible. $\Box$
\begin{corollary}\label{cor:cora} For $q\in T_n$, $p_n'(q)<1$.
\end{corollary}
\noindent Proof. The second derivative is positive in $T_n$, so the first derivative is always less
than its value of 1 at $n/(2n+1)$. $\Box$

While the original probability formulation makes no sense when $q=0$  (there is no optimal choice of $N$ when the underdog loses each play), the limit $\lim_{q\to 0}p_n(q)$, evident from Figure 3, is quite simple.
\begin{lemma}\label{lem:xx3} $\lim_{q\to 0^+}p_n(q)=\frac{1}{n+1}$
\end{lemma}
\noindent Proof.  Because $p_n'(q)>0$ (Lemma \ref{lem:xx2}), the values $p_n(q)$ decrease as $q\to 0$; the values are bounded and so the claimed limit exists. Thus we will use $p(0)$ to denote $\lim_{q\to 0^+}p_n(q)$.  Further, $0<p (0)<1$, for by Lemma \ref{lem:trap}, $p_n(q)>1/(2n+1)+q$ so
$p(0)\ge 1/(2n+1)$; and $p(0)<1-n/(2n+1)$ because the derivative is positive and $p_n(n/(2n+1))<1$ (Lemma \ref{lem:fact6}). Now, the derivative formula tells us that
\[\frac{((p-1) p (n p-n q+p-1))}{((q-1) q (n (q-p)+q))}-p'(q)=0,\]
 where $p$ denotes $p_n$. Multiplying both sides by the denominator turns this into
\[(p-1) p (n p-n q+p-1)-(q-1) q (n (q-p)+q) p'(q)=0.\]
 But the boundedness of $p_n'$ means that the limit of $q p'(q)$ is $0$, giving
\[\lim_{q\to 0^+}p(p-1)(np+p-1)=0,\]
 or
       $(p(0)-1) p(0) (p(0)(n+1)-1)=0$. Only the last factor can vanish, giving $p (0)=1/(n+1)$. $\Box$
\bgs
\begin{corollary}\label{cor:slp} (Slope convergence) The limit $\lim_{q\to 0^+}p_n'(q)$ exists.\end{corollary}
\noindent Proof. The slopes are bounded below (Lemma \ref{lem:xx2}) and monotonic by the convexity of $p_n$. $\Box$
\begin{corollary}\label{cor:slm} (Slope-limit formula)  $\lim_{q\to 0^+}p_n'(q)=n/(2(n+1)).$
\end{corollary}
\noindent Proof. Let $S$ denote the limit. Because the derivative formula $p'=((1-p) p (n p-n q+p-1))/((1-q) q (n  q-n p +q))$ holds for the slopes, we want the limit of this expression as $q\to 0$. Now, as $q\to 0$, $p\to 1/(n+1)$. So we can look at the numerator and denominator separately and see that we have a $0/0$ form, and we can use l'Hopital's rule to get the limit, where we use $p (q)$ for $p$. Forming the l'Hopital quotient, using $\lim p(q) = 1/(n+1)$, $\lim p'(q)=S$, and then letting $q$ be 0 yields $(n-(n-1)S)/(n+1)$. Setting this to equal to $S$
and solving gives the formula. $\Box$
\bgs
Corollary \ref{cor:slm} allows us to think of $p_n$ as a $C^1$ function on all of ${\mathbb R}$ as follows. Let $S_n$ be the limit of the slopes that the corollary provides. Define $p_n^*$ to agree with $p_n$ on $[0,n/(2n+1))$ and to be the linear function through $(0,p_n(0))$ of slope $S_n$ on $(-\infty,0]$, and the similar tangent-line extension on the right. It is easy to see using the Mean Value Theorem and the limit definition of $S_n$ that the limit of the slopes of the secants connecting $(0,p_n(0))$ to $(q,p_n(q))$ is $S_n$, giving the continuous differentiability of $p_n^*$.

\begin{corollary}\label{cor:slb} (Slope bound) $p_n'(q)>\frac{n}{2(n+1)}$ in $T_n$.\end{corollary}
\noindent Proof. By the previous corollary, because the slopes are monotonically increasing, by convexity. $\Box$
\begin{corollary}\label{cor:llb}  (Linear lower bound) In $T_n$, $p_n>K_n(q)\eqdef \frac{1}{n+1}+\frac{n}{2(n+1)}q$, and
\[N(q,p)\ge \lceil{(1-p)/(p-q/2)}\rceil\].
\end{corollary}
\noindent Proof. The linear function $K_n$ agrees with $p_n$ at $q=0$ by Lemma \ref{lem:xx3}. If $p_n$ dipped below $K_n$, the Mean Value Theorem would provide a
point contradicting Corollary \ref{cor:slb}. Inverting the bound on $p$ gives a bound on $N$. $\Box$
\bgs
The upper bound on $N$ given by Theorem \ref{th:ubn} and the piecewise lower bound obtained by combining Theorem \ref{th:low} and Corollary \ref{cor:llb} are useful computationally. For when the two bounds agree we know immediately that $N(q,p)=\lfloor{1/(2 (p-q))+1/2}\rfloor$.
 And in cases where $N$ equals the piecewise lower bound, that can be verified by a single
  $J$-evaluation: just check that $J_n(q,p)$  is negative when $n$ is the lower bound.
  The subset of $T$ for which the two bounds agree is a union of triangles -- the green region
   in Figure 9 -- and the total area of these triangles can be determined by integrating
   Mathematica's \texttt{Boole} function to get an expression for each level and then summing
    the results symbolically. The total area of the region in which $N$ equals the lower bound
     can be approximated by experimentation using thousands of points. The results are summarized
      in the next theorem.
\begin{theorem}\label{th:area}
\begin{enumerate}
\item For $(q,p)$ chosen uniformly from $T$, the probability that the upper bound of Theorem
\ref{th:ubn} and the combined lower bounds of Theorem \ref{th:low} and Corollary \ref{cor:slb} agree is
\[
(1-i) \psi (3-i)+(1+i) \psi (3+i)+\frac{\pi ^2}{4}+2 \gamma -\frac{115}{27},\]
where $\psi$ is the digamma function $\Gamma'/\Gamma$, and $\gamma$ is Euler's constant. This is roughly 0.60.
\item For $(q,p)$ as above, the probability that
\begin{equation}\label{eq:nbd}
N(q,p)=\max \left(\left\lceil{
   \frac{1-p}{p-\frac{q}{2}}}\right\rceil ,\left\lfloor{\frac{1}{2(p-q)}-\frac{1}{2}}\right\rfloor \right)
\end{equation}
 is approximately 0.87.
\end{enumerate}
\end{theorem}
\begin{corollary}\label{cor:imp} (Improved bounds on N) In $T$, $N^-(q,p)\le N(q,p)\le N^+(q,p)$, where
$N^-$ and $N^+$ are, respectively, the positive-radical solutions, $n$, to the quadratic equations:
\[2 n^2 (p-q)^3+2 n \left(p^2-3 p q-p+q^2+2 q\right)(p-q)-(1-p) q (1-3 p+3 q)=0,\]
and the equation DF$=n/(2(n+1))$, where DF is the derivative formula of Lemma \ref{lem:fact3}.
\end{corollary}
\noindent Proof. For the lower bound, the equation is equivalent to setting the second derivative formula (Lemma \ref{lem:xx1})
 to $0$ and clearing nonzero terms. First define
 \[ker=2 n^2 (p-q)^3+2 n \left(p^3-p^2 (4 q+1)+p q (4 q+3)-q^2 (q+2)\right)-(p-1) q (3p-3q-1),\]
 the result of clearing nonzero terms from the second derivative formula. Then the second derivative formula vanishes iff $ker$ does. Further, $\partial_nker>0$ in $T$ when we add the condition $L_n(q)<p$. Because $N(q,p)$ is not less than the least $n$ such that $p_n''(q)=0$, the result follows. The upper bound is obtained the same way, using the fact that the slope is not less than $n/(2(n+1))$, which becomes a quadratic relation. $\Box$
\bgs
Expanding the rational expression for $N^-$ in a series in powers of $p-q$ shows that
 it equals $1/(2(p-q))-3/2+1/(4 p(1-p))+O(p-q)$, which relates it nicely to the simpler bound of
  Theorem \ref{th:low}. Define $H(q,p)=\lceil{1/(2(p-q))-3/2+1/(4 p(1-p))}\rceil$; while $H$ is not a lower bound on $N$, it is a useful approximation when $p$ is close to $q$ and appears to be asymptotically perfect when the domain is rescaled (see subsection \ref{subs:9pt1}).

One can view the quadratic equations that give $N^{\pm}$ as cubic equations in $p$, in which case solving gives radical expressions for $p_n^{\pm}$, which bracket the curve $p_n$. Such bounds are useful when generating images such as Figure \ref{fig:fig10},  and also theoretically, as they are used in the proof of Theorem \ref{th:xxx1}.

We can measure how good the improved bounds are by assuming that the point $(q,p)$ is
 uniformly distributed in $T$. Then one can ask: (1) How often does $N^–=N^+$?
 (2) How often does $N^–=N$? The answers are: ``remarkably often.'' Figure 10 shows
  points in $T$ colored green if both bounds agree, red if the lower bound is correct, and yellow
   otherwise. The upper bound is sometimes correct, but not often. Of course, when the bounds
   agree we know $N$ immediately, and when the lower bound is correct, then that can be proved
    by verifying that $J_{N^-}(q,p)<0$ (more on computing $J$ in section 10 below). The blue
     curves are the graphs of $p_1^{\pm}$; they bracket $p_1$, defined by the yellow-red boundary.
       In Figure 10 the green area — where the two bounds agree — is
 75\% of the triangle and the green-plus-red area — where $N$ equals the lower bound — is 97\% of
  the area.
\begin{figure}[htp]\label{fig:fig9}
\begin{center}
\includegraphics[width=3in]{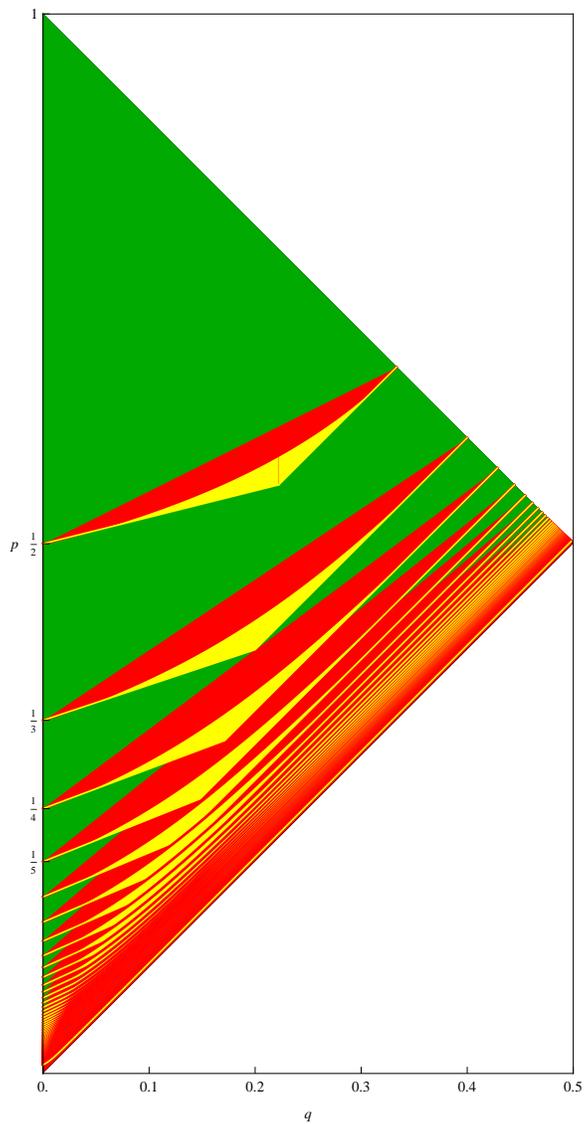}
\end{center}
\caption{The green region is where the two bounds on $N$ agree; the red region is where $N$ agrees with the lower bound. Thus in the combined red and green regions $N(q,p)$ is expressible exactly as $\max \left(\left\lceil\frac{1-p}{p-\frac{q}{2}}\right\rceil ,\left\lfloor\frac{1}{2(p-q)}-\frac{1}{2}\right\rfloor \right)$.}
\end{figure}
\begin{figure}[htp]\label{fig:fig10}
\begin{center}
\includegraphics[width=3in]{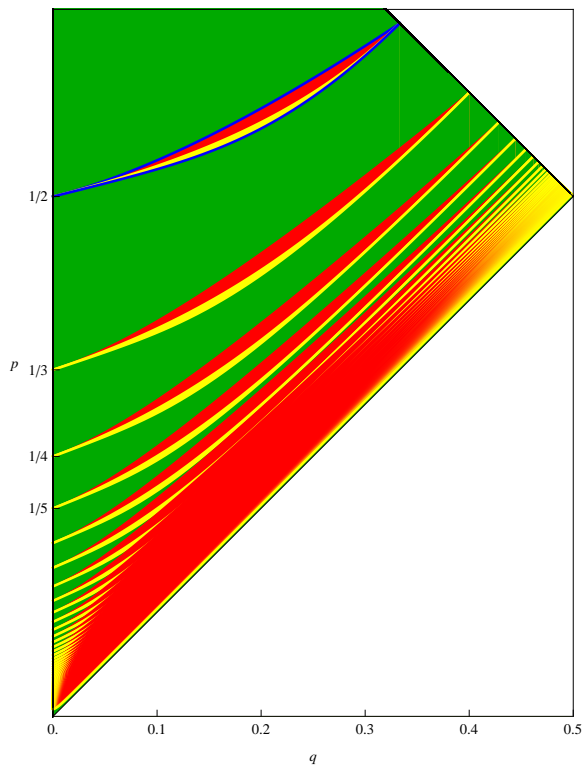}
\end{center}
\caption{The green region is where $N^-=N^+$; red is where $N^-=N<N^+$, and yellow is the rest. The blue curves are the graphs of $p_1^-$ and $p_1^+$.}
\end{figure}
\subsection{A harmonic rescaling}\label{subs:9pt1}
The views of Figures 9 and 10 do not show clearly what happens near the $p=q$ line. We can take a microscope to that area by harmonically rescaling the domain. We do this by first rotating $T$ clockwise $45^{\circ}$ and then stretching out the vertical scale: precisely, after the rotation we change each $y$-coordinate to $-(1/(2\sqrt{2} y))-1/2$. This is essentially just a change of coordinates from $(q,p)$ to $(p+q,1/(p-q))$; Figures 11 and 12 show the view through this microscope. The two approximations $N^-$ and $H$ exactly equal $N$ a large percentage of the time. For the region defined by $N\le 5000$ we found that $N=N^-$ in $>99\%$ of the region while $N$ and $H$ agree 96\% of the time. Thus we can conjecture that the probability of either equation holding is asymptotically 1.
\begin{figure}[htp]\label{fig:fig11}
\begin{center}
\includegraphics[width=3.5in]{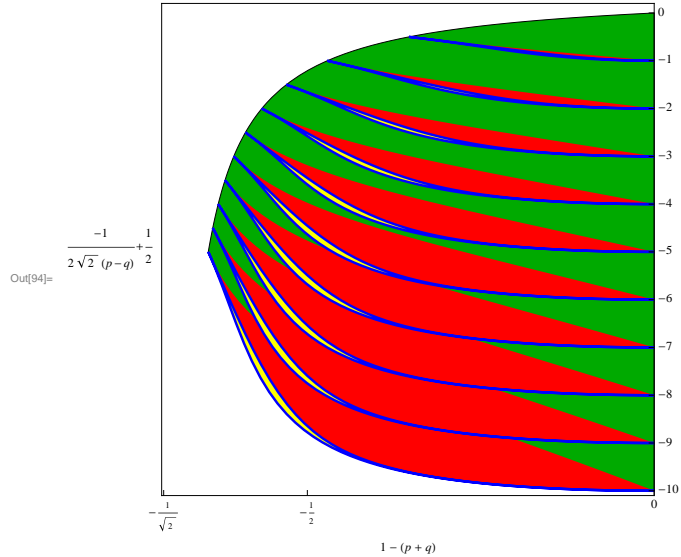}
\end{center}
\caption{ A rotated and vertically stretched view of the $(q,p)$ domain. The colors have the same meanings as in Figure 10. We see here that the region (yellow) where $N$ is not equal to $N^-$ is very small.}
\end{figure}
\begin{figure}[htp]\label{fig:fig12}
\begin{center}
\includegraphics[width=3.5in]{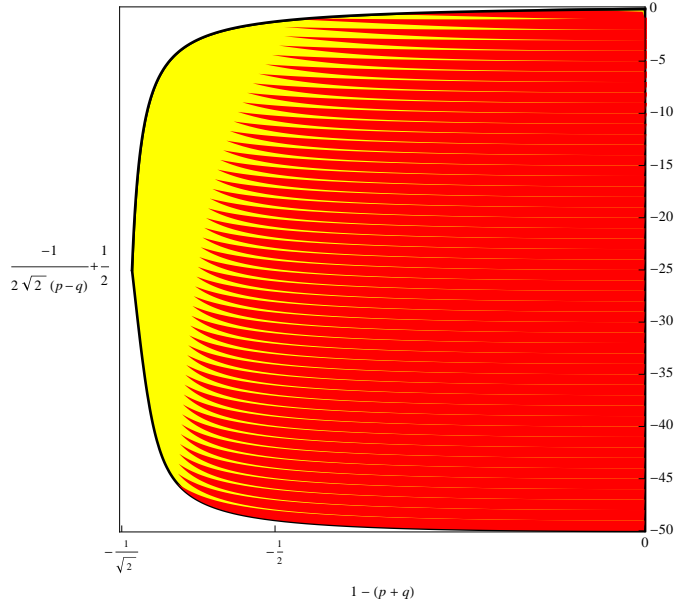}
\end{center}
\caption{A rotated and vertically stretched view of the $(q,p)$ domain in the region where $N\le 50$. Here red is where $N (q,p)=H (q,p)$, yellow where they are not equal. The proportion of the space in which $H$ is correct is 76\%, but it appears to converge to 100\% as $N\to\infty$.}
\end{figure}

\subsection{The situation when $p$ and $q$ are close}

The various diagrams suggest we look more closely at the situation near the $p=q$ line. The structure can be discerned by looking at the regions in which the excess of $N (q,p)$ over the simplest lower bound (Theorem \ref{th:bounds}) is constant. From this we will obtain the result that for any point $(q,q)$ there is an integer $\delta$ such that, close to $(q,q)$, we have $N(q,p)-\lfloor{\frac{1}{2(p-q)}+\frac12}\rfloor\le \delta$.

\medskip

\noindent\textbf{Definition.} For $(q,p)\in T$ let $\Delta(q,p)=N(q,p)-\lfloor{\frac{1}{2(p-q)}+\frac12}\rfloor$. the amount by which the optimal $n$ exceeds the lower bound derived from $L_n\le p_n$.

\medskip

Suppose $(q,p)$ is such that $p_{n+j}(q)<p<L_n(q)$ and also $L_{n+1}(q)<p$. The first inequality tells us that $n+1\le N\le n+j$. The last inequality means that the lower bound derived from $L_n$ is $n+1$. So $\Delta\le j-1$. Now we can profitably examine the regions determined by the intersection points of the $L$-lines with the nullclines $p_{n+j}$. Getting the intersection points is easy numerically using the function $p_n$, computed by a differential equation (section \ref{sec:sec10}); they are plotted as large black joined points in Figure 13. Note that each $L_n$ is tangent to $p_n$ at its right edge, strikes $p_{2n}$ at its left, and, because the slope of each nullcline is under 1 (Lemma \ref{lem:xx2}), hits each in-between $p$-graph in exactly one point.

Figure 13 tells the story. In the left image the colored regions correspond to constant values of $\Delta$. The uppermost black arc connects points common to: $L_1$ and $p_2$; $L_2$ and $p_3$; and in general $L_n$  and $p_{n+1}$. The second-highest black arc connects points common to: $L_2$ and $p_4$; $L_3$ and $p_5$; $L_4$ and $p_6$; in general $L_{n+1}$ and $p_{n+3}$. These arcs, which we cannot compute without computing $p_n$, divide $T$ into regions (right image of Figure 13) in each of which $\Delta$ takes on only two values.

\begin{figure}[]\label{fig:fig13}
\begin{center}
\includegraphics[width=6.5in]{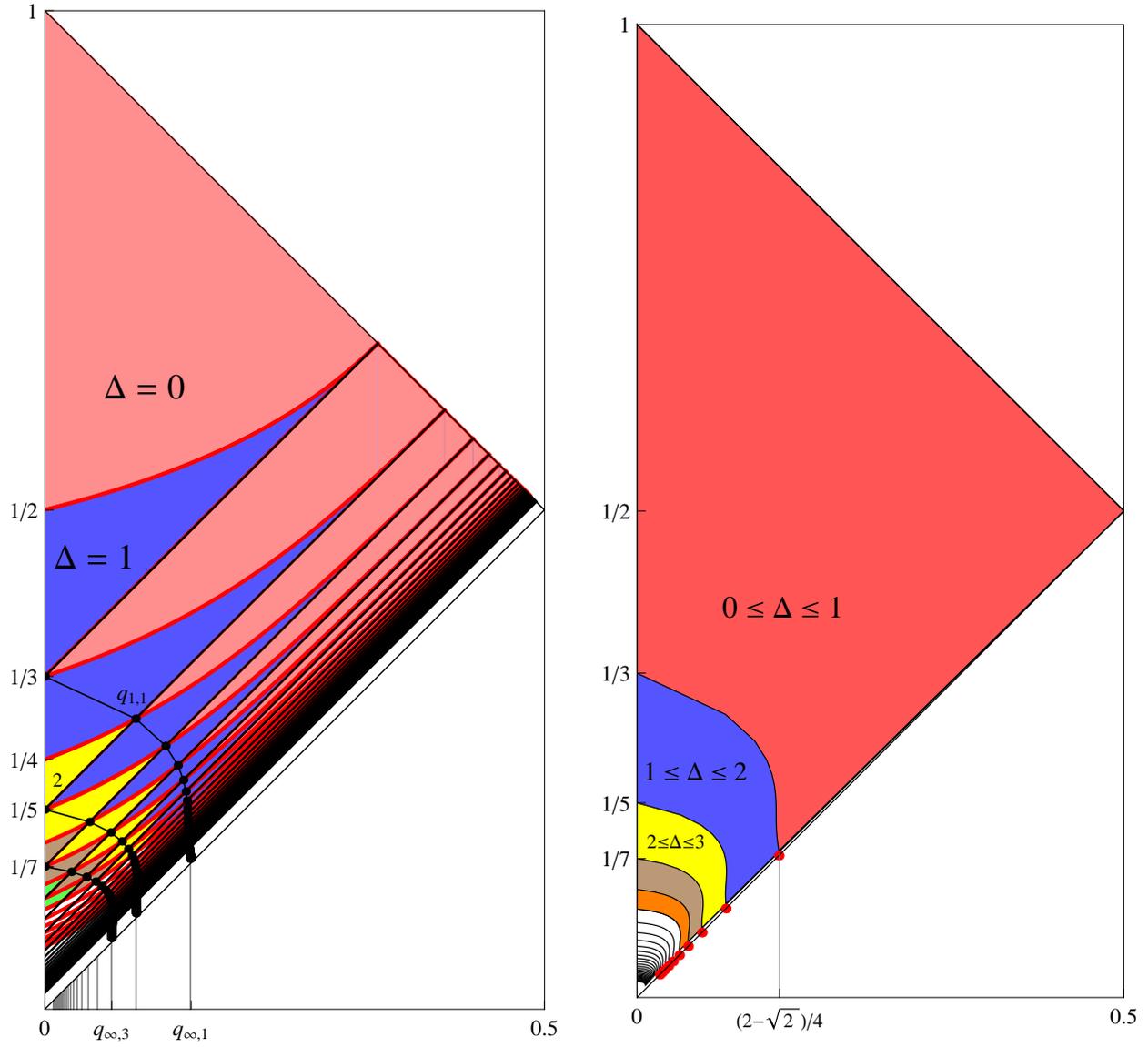}
\end{center}
\caption{The left image shows regions where $\Delta$ is constant (between the red $J$-nullclines and black lines $L_n$). The right image shows how the dividing arcs separate $T$ into regions where the discrepancy from the lower bound takes on one of two values.}
\end{figure}

\medskip

\noindent\textbf{Definition}.  For $n,j\ge 1$, let $q_{n,j}$ be the $q$-value of the intersection point of $L_{n+j+1}$ and $p_{n+2j-1}$ and let $q_{n,j}^-$ be the same with with $p$ replaced by the lower bound $p_{n+2j-1}^-$ derived from the convexity theorem (Corollary \ref{cor:imp}).

\medskip

Thus the arcs in Figure 13 are obtained by fixing $j$ and joining the points determined by $\{q_{n,j}\}$ as $n=1,2,3,\dots$. Because we used a lower bound on $p$, we know that $q_{n,j}^-<q_{n,j}$. Further, we can express $q_{n,j}^-$ quite simply by solving the equation $ker=0$ (see Corollary \ref{cor:imp}) after substituting $p=L_{n+j-1}(q)$,
\[q_{n,j}=\frac{-\sqrt{\frac{j \left(j^2+j (2 n-1)+(n-1)^2\right)}{j+1}}+j+n-1}{2
   j+2 n-1}.\]
The limit of this expression is easy to find; we use it to define $q_{\infty,j}=\frac12(1-\sqrt{j/(j+1)})$. A derivative shows easily that $q_{n,j}$ is increasing in $n$, and therefore $q_{n,j}$ approaches $q_{\infty,j}$ from the left. The points $(q_{\infty,j},q_{\infty,j})$ are the limits of the arcs in Figure 13; we now turn to a proof of this useful fact. Because $q_{n,j}^- <q_{n,j}$ and approaches the limit from the left, we need only show that $q_{n,j}<q_{\infty,j}$. This inequality is not true in general (though it does appear to be true when $j=1$) but we can show that, for each fixed $j$, it is true for sufficiently large $n$, and that suffices for the limit.

\begin{theorem}\label{th:xxx1}
 For fixed $j$ and sufficiently large $n$, $q_{n,j}<q_{\infty,j}$.\end{theorem}
\noindent Proof. To say that $q_{n,j}$ is to the left of $q_{\infty,j}$ is the same as saying that $L_{n+j-1}(q_{\infty,j})$ lies above $p_{n+2j-1}(q_{\infty,j})$. This in turn is the same as the assertion that $J_{n+2j-1}(q_{\infty,j},L_{n+j-1}(q_{\infty,j}))<0$. And this is equivalent to $r_{n+2j}(u)>\rho$, where $q_{\infty,j}$ and $L_{n+j-1}(q_{\infty,j})$ are used for $q$ and $p$ in $\rho$ and $u$.

Now define $g(n)={2n\choose n}4^{-n}$ and let
\[ P_{n,m}^*(z)=z^ng(n)\sum_{\nu=0}^m g(\nu)\left(\prod_{j=1}^{\nu}\frac{2j-1}{2n-1j+1}\right)z^{-2\nu}(1-z^{-2})^{-\nu-\frac12}.
\]
This is the asymptotic series of Laplace-Heine \cite[Thm. 8.21.3]{gs} for the Legendre polynomials $P_n(u)$ for $u>1$. Here $z=u+\sqrt{u^2-1}$. Thus
\[P_n(u)=P_{n,m}^*(u)+O(n^{-m-\frac32}z^n).\]
Therefore, forming the ratio, we have
\[r_n(u)=\frac{P_{n,m}^*(u)}{P_{n-1,m}^*(u)}+O(n^{-m-\frac32}).\]
To show that $r_{n+2j}(u)-\rho>0$ we must show that
\[r_{n+2j}\left(\frac{2n+4jn-1+4j^2-2\sqrt{j(j+1)}}{2(1+j+\sqrt{j(j+1)}(2n+2j-1))}\right)>
\frac{2(j+1)(n+j-1)}{\sqrt{j(j+1)}(2n+2j-1)-(j+1)}.\]
We can now apply the asymptotic series using three terms ($m=3$), and take three terms of the Taylor series of the result, centered at $\infty$. Fewer than three terms are not sufficient. When we do this using \textit{Mathematica}'s \texttt{Series} command and some further simplification we get the following expression, whose positivity concludes the proof-
\[r_{n+2j}(u)-\rho=\frac12(j+1)^2\left(1-\sqrt{j/(j+1)}\right)n^{-3}+O(n^{-4}).\ \Box .\]

Given $q$; let $j_0(q)=\lceil{(1-2q)^2/(4q(1-q))}\rceil$, the largest $j$ such that $q_{\infty,j}\le q$, let $n_0(j)$ be the least integer $n$ guaranteed by the theorem, i.e., $q_{n,j}<q_{\infty,j}$ for $n\ge n_0$. It appears that $n_0<3j^2$, but we have no proved bound.
\begin{corollary}\label{cor:newub}
If $q<p<L_{n_0(j_0(q))}$ then $N(q,p)\le \lfloor{\frac{1}{2(p-q)}+\frac12}\rfloor+\lceil{\frac{(1-2q)^2}{4q(1-q)}}\rceil $.
\end{corollary}
\noindent Proof. The definition of $n_0$ means that as one starts at the point $(q,L_{n_0(j_0(q))})$ and moves down, one strikes, in alternating order, the graphs of $p_n$ and the lines $L_n$. This means that for these points $\Delta$  equals its value just under $(q,L_{n_0(j_0(q))})$, or 1 greater. Since the value of $\Delta$ just below the line $L_{n_0}$ is at most $\lceil{(1-2 q)^2/(4 (1-q) q)}\rceil-1$ (see comments at start of this section), the result is proved. $\Box$
\begin{corollary}
Given $q_0$, with $0<q_0\le \frac12$, we have
$\lim_{(p,q)\to (q_0,q_0)}(q-p)N(q,p)=\frac12$
\end{corollary}
\noindent Proof. Immediate from the previous corollary, which implies that
\[\frac{1}{(2 (p-q))}-\frac12\le N(q,p)\le \frac{1}{2 (p-q)}+\frac{(1-2 q)^2}{4 (1-q) q}+2.\ \Box \]

If in Theorem \ref{th:xxx1}  one replaces the sharp $q_{\infty,j}$ by simply $1/(2j+1)$ one obtains a much weaker theorem, which can be proved in the same way; i.e., $q_{n,j}<1/(2j+1)$ for sufficiently large $n$. However, in this case it appears that the result is true for all $n$, as is evident from Figure 13. A proof of this would yield a new upper bound on $N (q,p)$, weaker than the one in Corollary \ref{cor:newub}, but with the advantage of being true for all $(q,p)$.

\medskip

\noindent \textbf{Open problems}. Prove that $r_{n+2j}(u)<\rho$, where $1/(2j+1)$ and $L_{n+j-1}(1/(2j+1))$ are, respectively, used for $q$ and $p$ in $\rho$ and $u$. Find estimates for $n_0(j)$. Show that $n_0(1)=1$.

\section{An algorithm for computing $N (q,p)$}\label{sec:sec10}
A straightforward algorithm to compute $N (q,p)$ first uses symmetry to restrict to $T$ and then finds the smallest integer $n$ such that $p\ge p_n(q)$; that value of $n$ is $N (q,p)$ by Lemma \ref{lem:lem13}. One can start with the simple or improved bounds and then use either bisection or the secant method,  repeatedly checking whether $J_n(q,p)$ is positive or negative. This method works fine when $N$ is of modest size, but when $N$ is large the Legendre polynomials cannot be explicitly computed. A solution is to use the integral formula given in (\ref{eq:intrep}), which is a fine substitute for $P_n$. That formula means that we can determine the sign of $J_n$ for each trial by using numerical integration on
\[(1-p+q)\int_0^{\pi}\left(u+\sqrt{u^2-1}\cos{t}\right)^ndt-
(1-p-q)\int_0^{\pi}\left(u+\sqrt{u^2-1}\cos{t}\right)^{n+1}dt.\]

Of course, high-precision must be used as appropriate. One needs enough accuracy to account for the full precision of $n$ which will be used as a trial in the root-finding process. Further, the expression $\sqrt{u^2-1}$ can be numerically unstable for extreme values of $p$ and $q$ and one should use the equivalent form $2\sqrt{(1-p) p (1-q) q}/(1-p-q)$.

But one needs only the sign of the integral above. In Mathematica this means that when computing the integral numerically one needs a large working precision, but the accuracy goal can be quite small. The method is  robust and takes only a few seconds to compute $N(10^{-100}, 2\cdot 10^{-100})$,
which is

\[72768\,90317\,94675\,98852\,95987\,53552\,38752\,84521\,10838\,88022\,00705\,
28794\,63897\,19626\,49789\]
\nobreak
\vspace{-.2in}
\nobreak
\noindent\,77512\,24788\,32188\,39061\,36928.\\

{\renewcommand{\arraystretch}{1.5}
\begin{center}
  \begin{tabular}{| c | c | c | }
    \hline
    $q$&$p$&$N(q,p)$\\ \hline\hline
    $10^{-5}$ & $2\cdot 10^{-5}$ & 72768\\ \hline
    $10^{-10}$ & $2\cdot 10^{-10}$ & 7276890317\\ \hline
    $10^{-15}$ & $2\cdot 10^{-15}$ & 727689031794675\\ \hline
    $10^{-20}$ & $2\cdot 10^{-20}$ & 72768903179467598852\\ \hline
    $10^{-25}$ & $2\cdot 10^{-25}$ & 7276890317946759885295987\\ \hline
    $10^{-30}$ & $2\cdot 10^{-30}$ & 727689031794675988529598753552\\ \hline
  \end{tabular}
\bigskip

  Table 1. The integration algorithm allows one to get giant values of $N (q,p)$.
\end{center}}

When one wants not just the sign of $J_n$, but a numerical approximation to the
full $J_n$-nullcline -- the graph of $p_n$ -- one can use a numerical differential
 equation approach. Because of the derivative formula and the known values at $0$, we can set up
  the initial-value problem as $p(0) = 1/(n+1)$ and
  \[ p'(q)=\frac{p(1-p)(np-nq+p-1)}{q(1-q)(n(q-p)+q)}\]
 if $q>0$ and $n/(2(n+1))$ if $q=0$. This approach is a quick way to generate graphs of
 $p_n$, such as those shown in various figures in this paper. It can also be used in an algorithm
  for computing $N (q,p)$ where it can sometimes be faster than the use of numerical integration
  because the solution is needed only on the interval $[0,q]$, while the integrals are computed
   from $0$ to $\pi$.

\section{Some open questions}
\begin{enumerate}
\item Improve the bounds on $N (q,p)$.
\item  Find a more efficient algorithm for computing $N(q,p)$ when $(q,p)$ is near the origin.
\item Generalize, in a natural way, these results to the case of three players.
\item Prove the second derivative conjectures, which we obtained heuristically by manipulating Taylor polynomials:
\begin{enumerate}
\item \[\lim_{q\to 0^+}p_n''(q)=\frac{2n^2+5n+2}{6(n+1)},\]
 and
 \item \[\lim_{q\to \left(\frac{n}{2n+1}\right)^-}p_n''(q)=\frac{4(2n+1)}{n(2n^2+n-1)}.\]
 \end{enumerate}
\item (Asymptotic closed form conjectures)
\begin{enumerate}
\item If $q$ and $p$ are chosen uniformly from the harmonically rescaled domain then the probability that $N (q,p)=N^-(q,p)$ approaches 1 as $N\to\infty$.
\item Same as above, but with the conclusion that $N=\lceil{1/(2(p-q))-3/2+1/(4 p(1-p))}\rceil$ with asymptotic probability 1.
\end{enumerate}
\end{enumerate}

\vspace{.5 in}

\noindent Macalester College, St. Paul, MN, 55105;\ \texttt{<addona@macalester.edu>}\\
\noindent Macalester College, St. Paul, MN, 55105;\ \texttt{<wagon@macalester.edu>}\\
\noindent University of Pennsylvania, Philadelphia, PA 19104-6395;\ \texttt{<wilf@math.upenn.edu>}

\end{document}